\documentclass[journal,twoside]{IEEEtran}
\usepackage{cite}
\ifCLASSINFOpdf
\else
\fi
\usepackage[cmex10]{amsmath}
\usepackage{array}
\usepackage{mdwmath}
\usepackage{mdwtab}
\usepackage[caption=false,font=footnotesize]{subfig}
\usepackage{stfloats}
\usepackage{url}
\usepackage[american,fulldiode]{circuitikz} %%Circuit drawing tools
\usepackage{graphicx} %%For loading graphic files
\usepackage{multirow}
\usepackage{amsfonts}
\usepackage{bbm}
\usepackage{algorithm}
\usepackage{algpseudocode}
\usepackage{multirow}
\usepackage[section]{placeins}
\usepackage{epstopdf}

\usepackage{slashbox}
\usepackage{amssymb}

\usepackage{lipsum,amsmath,multicol}

\renewcommand\Re{\operatorname{Re\mathfrak{}}}
\renewcommand\Im{\operatorname{Im\mathfrak{}}}

\newcommand{\matpower}{M{\sc atpower}}
\usepackage{multirow}

\usepackage{float}
\usepackage{tablefootnote}
\usepackage[multiple]{footmisc}

\makeatletter
\def\ps@headings{%
\def\@oddhead{\mbox{}\scriptsize\rightmark \hfil \thepage}%
\def\@evenhead{\scriptsize\thepage \hfil \leftmark\mbox{}}%
\def\@oddfoot{}%
\def\@evenfoot{}}
\def\BState{\State\hskip-\ALG@thistlm}
\makeatother

\renewcommand\Re{\operatorname{Re\mathfrak{}}}
\renewcommand\Im{\operatorname{Im\mathfrak{}}}

\allowdisplaybreaks

\usepackage{varwidth}
\newcolumntype{M}[1]{>{\begin{varwidth}[t]{#1}}l<{\end{varwidth}}}

\newcommand\NoDo{\renewcommand\algorithmicdo{}}

\begin{document}
\title{A Laplacian-Based Approach for Finding Near Globally Optimal Solutions to OPF Problems}
\author{Daniel K. Molzahn, \IEEEmembership{Member, IEEE}, C\'{e}dric Josz, Ian A. Hiskens, \IEEEmembership{Fellow, IEEE}, \\ and Patrick Panciatici \IEEEmembership{Member, IEEE}% <-this % stops a space
\thanks{Argonne National Laboratory, Energy Systems Division: dmolzahn@anl.gov,
University of Michigan, Dept. of Electrical Engineering and Computer Science: hiskens@umich.edu,
R\'eseau de Transport d'\'Electricit\'e, R\&D Dept., cedric.josz@rte-france.com, patrick.panciatici@rte-france.com}}% <-this % stops a space
\maketitle

\begin{abstract}
A semidefinite programming (SDP) relaxation globally solves many optimal power flow (OPF) problems. For other OPF problems where the SDP relaxation only provides a lower bound on the objective value rather than the globally optimal decision variables, recent literature has proposed a penalization approach to find feasible points that are often nearly globally optimal. A disadvantage of this penalization approach is the need to specify penalty parameters. This paper presents an alternative approach that algorithmically determines a penalization appropriate for many OPF problems. The proposed approach constrains the generation cost to be close to the lower bound from the SDP relaxation. The objective function is specified using iteratively determined weights for a Laplacian matrix. This approach yields feasible points to the OPF problem that are guaranteed to have objective values near the global optimum due to the constraint on generation cost. The proposed approach is demonstrated on both small OPF problems and a variety of large test cases representing portions of European power systems.
\end{abstract}



% \fi
\IEEEpeerreviewmaketitle

\section{Introduction}
\label{l:introduction}

\IEEEPARstart{T}{he} optimal power flow (OPF) problem determines an
optimal operating point for an electric power system in terms of a
specified objective function (typically generation cost per unit
time). Equality constraints for the OPF problem are dictated by the
network physics (i.e., the power flow equations) and inequality
constraints are determined by engineering limits (e.g., voltage
magnitudes, line flows, and generator outputs).

The OPF problem is non-convex due to the non-linear power flow
equations, may have local optima~\cite{bukhsh_tps}, and is generally NP-Hard~\cite{lavaei_tps}, even for relatively simple cases such as tree-topologies~\cite{NPhard}. There is a large literature on solving OPF problems using local optimization techniques (e.g., successive quadratic programs, Lagrangian relaxation, heuristic optimization, and interior point methods~\cite{opf_litreview1993IandII,ferc4}).  These techniques are generally well suited to solving large problems. However, while local solution techniques often find global solutions~\cite{molzahn_lesieutre_demarco-global_optimality_condition}, they may fail to converge or converge to a local optimum~\cite{bukhsh_tps,ferc5}. Furthermore, they are unable to quantify solution optimality relative to the global solution.

There has been significant recent research focused on convex
relaxations of OPF problems. Convex relaxations lower bound the
optimal objective value and can certify infeasibility of OPF
problems. These relaxations also yield the globally optimal decision
variables for many OPF problems (i.e., the relaxations are often
``exact''). Second-order cone programming (SOCP) relaxations can globally solve
OPF problems for radial networks that satisfy certain non-trivial
technical conditions~\cite{low_tutorial}. Semidefinite programming
(SDP) relaxations globally solve a broader class of OPF
problems~\cite{bai2008,lavaei_tps,molzahn_holzer_lesieutre_demarco-large_scale_sdp_opf}.

Recently, the SDP relaxation has been generalized to a family of ``moment'' relaxations using the Lasserre hierarchy for polynomial optimization~\cite{pscc2014,cedric_msdp,ibm_paper}. With increasing relaxation order, the moment relaxations globally solve a broader class of OPF problems at the computational cost of larger SDPs.

By exploiting network sparsity and selectively applying the computationally intensive higher-order constraints, the moment relaxations are capable of globally solving larger OPF problems~\cite{molzahn_hiskens-sparse_moment_opf}, including problems with several thousand buses representing portions of European power systems~\cite{cdc2015}. Ongoing efforts include further increasing computational speed. Recent work includes implementing the higher-order moment constraints with a faster SOCP formulation~\cite{powertech2015} and development of a complex version of the Lasserre hierarchy~\cite{josz_molzahn-complex_hierarchy}.

As an alternative to the moment relaxations, other literature has
proposed an objective function penalization approach for finding
feasible points that are near the global optimum for the OPF
problem~\cite{lavaei_mesh,lavaei_allerton2014}. The penalization
approach has the advantage of not using potentially
computationally expensive higher-order moment constraints, but has the
disadvantage of requiring the choice of appropriate penalization
parameters. This choice involves a compromise, as the
  parameters must induce a feasible solution to the original problem
  while avoiding large modifications to the problem that would cause
  unacceptable deviation from the global optimum.

%The penalty parameters must be chosen to balance obtaining
%a feasible solution to the original problem while avoiding excessive
%modifications to the problem in order to stay close to the global
%optimum.

The penalization formulation in the existing
literature~\cite{lavaei_allerton2014} generally requires specifying
penalty parameters for both the total reactive power injection and
apparent power flows on certain lines. Penalty parameters
  in the literature range over several orders of magnitude for various
  test cases, and existing literature largely lacks  systematic algorithms for determining appropriate parameter values. Recent work~\cite{cdc2015} proposes a ``moment+penalization'' approach that eliminates the need to choose apparent power flow penalization parameters, but still requires selection of a penalty parameter associated with the total reactive power injection.

%This paper presents an approach that does not require selection of any penalty parameters and does not employ constraints from higher-order moment relaxations.

This paper presents an iterative algorithm that builds an objective function intended to yield near-globally-optimal solutions to OPF problems. The algorithm is applicable for cases where the SDP relaxation is not exact but has a small ``relaxation gap'' (i.e., the solution to the SDP relaxation has an objective value that is close to the true globally optimal objective value). The proposed algorithm first solves the SDP relaxation to obtain a lower bound on the optimal objective value. For many practical OPF problems, this lower bound is often very close to the global optimum. The proposed approach modifies the SDP relaxation by adding a constraint that the generation cost must be within a small percentage (e.g., 0.5\%) of this lower bound. This percentage is the single externally specified parameter in the proposed approach.

This constraint on the generation cost provides freedom to specify an
objective function that aims to obtain a \emph{feasible} rather than
\emph{minimum-cost} solution for the OPF problem. In other words, we
desire an objective function such that the SDP relaxation yields a
feasible solution to the original non-convex OPF problem, with
near-global optimality ensured by the constraint on generation
cost.

This paper proposes an algorithm for calculating an appropriate
objective function defined using a weighted Laplacian matrix. The
weights are determined iteratively based on the mismatch between
 the solution to the relaxation and the power flows resulting from a related set of voltages. The paper will formalize these concepts and
  demonstrate that this approach results in near global solutions to
many OPF problems, including large test cases. Like many penalization/regularization techniques~\cite{lavaei_mesh,lavaei_allerton2014}, the proposed approach is not guaranteed to yield a feasible solution.\footnote{The moment+penalization method in~\cite{cdc2015} applies the Lasserre hierarchy for polynomial optimization to OPF problems that are augmented with a reactive power penalty. For sufficiently high relaxation orders, the approach in~\cite{cdc2015} is guaranteed to yield a feasible solution. However, solving relaxations from the Lasserre hierarchy can be computationally challenging.} 
As supported by the results for several large-scale, realistic test cases, the proposed algorithm broadens the applicability of the SDP relaxation to achieve operating points for many OPF problems that are within specified tolerances for both constraint feasibility and global optimality.

There is related work that chooses the objective function of a
relaxation for the purpose of obtaining a feasible solution for the
original non-convex problem. For instance,~\cite{allerton2011} specifies objective functions that are
linear combinations of squared voltage magnitudes in order to find
multiple solutions to the power flow equations. Additionally,~\cite{lavaei_pf} proposes a method for determining an
objective function that yields solutions to the power flow equations
for a variety of parameter choices. The objective function
in~\cite{lavaei_pf} is defined by a matrix with three properties:
positive semidefiniteness, a simple eigenvalue of 0, and null space
containing the all-ones vector. We note that the weighted Laplacian
objective function developed in this paper is a special case of an objective function that also has these three properties.

This paper is organized as follows. Section~\ref{l:opf_formulation} introduces the OPF formulation studied in this paper. Section~\ref{l:sdp_relaxation} reviews the SDP relaxation from previous literature. Section~\ref{l:lap_pen} describes the Laplacian objective function approach that is the main contribution of this paper. Section~\ref{l:results} demonstrates the effectiveness of the proposed approach through application to a variety of small OPF problems as well as several large test cases representing portions of European power systems. Section~\ref{l:conclusion} concludes the paper.

\section{Optimal Power Flow Problem}
\label{l:opf_formulation}

We first present an OPF formulation in terms of complex voltage
coordinates, active and reactive power injections, and apparent power
line flow limits. Consider an $n$-bus system with $n_l$ lines, where
$\mathcal{N} = \left\lbrace 1, \ldots, n \right\rbrace$ is the set of
buses, $\mathcal{G}$ is the set of generator buses, and $\mathcal{L}$
is the set of lines. The network admittance matrix is $\mathbf{Y} =
\mathbf{G} + \mathbf{j} \mathbf{B}$, where $\mathbf{j}$ denotes the
imaginary unit. Let $P_{Dk} + \mathbf{j} Q_{Dk}$ represent the active
and reactive load demand and $V_k = V_{dk} + \mathbf{j} V_{qk}$ the
voltage phasors at each bus~$k \in \mathcal{N}$. Superscripts ``max''
and ``min'' denote specified upper and lower limits. Buses without
generators have maximum and minimum generation set to zero. Let
$c_{k2}$, $c_{k1}$, and $c_{k0}$ denote the coefficients of a convex
quadratic cost function for each generator $k \in \mathcal{G}$.

The power flow equations describe the network physics:
\begin{subequations}
%\small
\label{opf_balance}
\begin{align}\nonumber
P_{Gk} = & V_{dk} \sum_{i=1}^n \left( \mathbf{G}_{ik} V_{di} - \mathbf{B}_{ik} V_{qi} \right) &  &  \\[-5pt] 
\label{opf_Pbalance}  & + V_{qk} \sum_{i=1}^n \left( \mathbf{B}_{ik}V_{di} + \mathbf{G}_{ik}V_{qi} \right) + P_{Dk} \\ 
\nonumber Q_{Gk} = & V_{dk} \sum_{i=1}^n \left( -\mathbf{B}_{ik}V_{di} - \mathbf{G}_{ik} V_{qi}\right) \\[-5pt]
\label{opf_Qbalance} & + V_{qk} \sum_{i=1}^n \left( \mathbf{G}_{ik} V_{di} - \mathbf{B}_{ik} V_{qi}\right) + Q_{Dk}.
\end{align}
\end{subequations}

We use a line model with an ideal transformer that has a specified turns ratio $\tau_{lm} e^{\mathbf{j}\theta_{lm}} \colon 1$ in series with a $\Pi$ circuit with series impedance $R_{lm} + \mathbf{j} X_{lm}$ (equivalent to an admittance of $g_{lm} + \mathbf{j} b_{lm} = \frac{1}{R_{lm} + \mathbf{j} X_{lm}}$) and total shunt susceptance $\mathbf{j} b_{sh,lm}$. The line flow equations are
\begin{subequations}
\small
\begin{align}
\nonumber  & P_{lm} = \left( V_{dl}^2 + V_{ql}^2\right) g_{lm}/\tau_{lm}^2 \\ \nonumber
& \quad +
\left(V_{dl}V_{dm} + V_{ql}V_{qm}\right) \left(b_{lm}\sin\left(\theta_{lm} \right) - g_{lm}\cos\left(\theta_{lm} \right) \right) / \tau_{lm} \\
\label{Plm}& \quad +
\left(V_{dl}V_{qm} - V_{ql}V_{dm}\right) \left(g_{lm}\sin\left(\theta_{lm}\right) + b_{lm}\cos\left(\theta_{lm}\right)\right) / \tau_{lm}  \\ % \end{align} \begin{align}
\nonumber & P_{ml} = \left(V_{dm}^2 + V_{qm}^2 \right)g_{lm} \\ \nonumber & \quad -
\left(V_{dl}V_{dm} + V_{ql}V_{qm} \right)\left(g_{lm}\cos\left(\theta_{lm}\right) + b_{lm}\sin\left(\theta_{lm}\right) \right) / \tau_{lm} \\ \label{Pml}
& \quad + 
\left(V_{dl}V_{qm} - V_{ql}V_{dm} \right) \left(g_{lm}\sin\left(\theta_{lm}\right) - b_{lm}\cos\left(\theta_{lm}\right) \right) / \tau_{lm}\\
\nonumber & Q_{lm} = -\left( V_{dl}^2 + V_{ql}^2\right) \left(b_{lm} + \frac{b_{sh,lm}}{2}\right) / \tau_{lm}^2 \\ \nonumber 
& \quad +
\left(V_{dl}V_{dm} + V_{ql}V_{qm}\right) \left(b_{lm}\cos\left(\theta_{lm} \right) + g_{lm}\sin\left(\theta_{lm} \right) \right) / \tau_{lm} \\ \label{Qlm} & \quad +
\left(V_{dl}V_{qm} - V_{ql}V_{dm}\right) \left(g_{lm}\cos\left(\theta_{lm}\right) - b_{lm}\sin\left(\theta_{lm}\right)\right) / \tau_{lm} \\
\nonumber & Q_{ml} = -\left( V_{dm}^2 + V_{qm}^2\right) \left(b_{lm} + \frac{b_{sh,lm}}{2}\right) \\ \nonumber & \quad +
\left(V_{dl}V_{dm} + V_{ql}V_{qm}\right) \left(b_{lm}\cos\left(\theta_{lm} \right) - g_{lm}\sin\left(\theta_{lm} \right) \right)  / \tau_{lm} \\ \label{Qml} & \quad + 
\left(-V_{dl}V_{qm} + V_{ql}V_{dm}\right) \left(g_{lm}\cos\left(\theta_{lm}\right) + b_{lm}\sin\left(\theta_{lm}\right)\right) / \tau_{lm}. %\\
%\label{Slm} & \left(S_{lm}\right)^2 = \left(P_{lm}\right)^2 + \left(Q_{lm}\right)^2 \\
%\label{Sml} & \left(S_{ml}\right)^2 = \left(P_{ml}\right)^2 + \left(Q_{ml}\right)^2
\end{align}
\end{subequations}

The classical OPF problem is then
\begin{subequations}
\label{opf}
%\small
\begin{align}
\label{opf_obj} & \min_{V_d,V_q}\quad \sum_{k \in \mathcal{G}} c_{k2} P_{Gk}^2 + c_{k1} P_{Gk} + c_{k0} \qquad \mathrm{subject\; to} \hspace{-160pt} & \\
\label{opf_P} &  \quad P_{Gk}^{\mathrm{min}} \leq P_{Gk} \leq P_{Gk}^{\mathrm{max}} & \forall k \in \mathcal{N} \\
\label{opf_Q} &  \quad Q_{Gk}^{\mathrm{min}} \leq Q_{Gk} \leq Q_{Gk}^{\mathrm{max}} &  \forall k \in \mathcal{N} \\
\label{opf_V} &  \quad (V_{k}^{\mathrm{min}})^2 \leq V_{dk}^2 + V_{qk}^2 \leq (V_{k}^{\mathrm{max}})^2 &  \forall k \in \mathcal{N}  \\
\label{opf_Slm} & \quad \left(P_{lm}\right)^2 + \left(Q_{lm}\right)^2 \leq \left(S_{lm}^{\mathrm{max}}\right)^2 &  \forall \left(l,m\right) \in \mathcal{L} \\ 
\label{opf_Sml} & \quad \left(P_{ml}\right)^2 + \left(Q_{ml}\right)^2 \leq \left(S_{lm}^{\mathrm{max}}\right)^2 &  \forall \left(l,m\right) \in \mathcal{L} \\ 
\label{opf_Vref} & \quad V_{q1} = 0.
\end{align}
\end{subequations}
Constraint~\eqref{opf_Vref} sets the reference bus angle to zero.

\section{Semidefinite Relaxation of the OPF Problem}
\label{l:sdp_relaxation}

This section describes an SDP relaxation of the OPF problem adopted from~\cite{lavaei_tps,molzahn_holzer_lesieutre_demarco-large_scale_sdp_opf,andersen2014}. Let $e_k$ denote the
$k^{th}$ standard basis vector in $\mathbb{R}^n$. Define 
  $Y_k = e_k e_k^\intercal \mathbf{Y}$, where $\left(\cdot\right)^\intercal$ indicates the transpose operator.

Matrices employed in the bus power injection, voltage magnitude, and angle reference constraints are
\begin{subequations}
\begin{align}\label{Yi}\mathbf{Y}_k & = \frac{1}{2} \begin{bmatrix}
\Re \left(Y_k + Y_k^\intercal\right) & \Im \left(Y_k^\intercal - Y_k\right) \\
\Im \left(Y_k - Y_k^\intercal\right) & \Re \left(Y_k\ + Y_k^\intercal\right)
\end{bmatrix} \\
\label{Yibar}\mathbf{\bar{Y}}_k & = -\frac{1}{2} \begin{bmatrix}
\Im \left(Y_k + Y_k^\intercal\right) & \Re \left(Y_k - Y_k^\intercal\right) \\
\Re \left(Y_k^\intercal - Y_k\right) & \Im \left(Y_k + Y_k^\intercal\right)
\end{bmatrix} \\
\label{Mi}\mathbf{M}_k & = \begin{bmatrix}
e_k e_k^\intercal & \mathbf{0} \\
\mathbf{0} & e_k e_k^\intercal
\end{bmatrix} \\
\label{Ni}\mathbf{N}_k & = \begin{bmatrix}
\mathbf{0} & \mathbf{0} \\
\mathbf{0} & e_k e_k^\intercal
\end{bmatrix}
\end{align}
\end{subequations}
where $\Re\left(\cdot\right)$ and $\Im\left(\cdot\right)$ return the real and imaginary parts, respectively, of a complex argument.

Define $f_i$ as the $i^{th}$ standard basis vector in
$\mathbb{R}^{2n}$, and define:
\begin{subequations}
\begin{align} c_{lm} & = \left(g_{lm} \cos\left(\theta_{lm}\right)-b_{lm}\sin\left(\theta_{lm}\right)\right) / \left(2 \tau_{lm}\right)\\
c_{ml} & = \left(g_{lm} \cos\left(\theta_{lm}\right)+b_{lm}\sin\left(\theta_{lm}\right)\right) / \left(2 \tau_{lm}\right) \\
s_{lm} & = \left(g_{lm} \sin\left(\theta_{lm}\right)+b_{lm}\cos\left(\theta_{lm}\right)\right) / \left( 2 \tau_{lm}\right) \\
s_{ml} & = \left(g_{lm} \sin\left(\theta_{lm}\right)-b_{lm}\cos\left(\theta_{lm}\right)\right) / \left( 2\tau_{lm}\right).
\end{align}
\end{subequations}
Matrices employed in the line flow constraints can then be written:
\begin{subequations}
\begin{align} \nonumber \mathbf{Z}_{lm} & = \frac{g_{lm}}{\tau_{lm}^2} \left(f_l^{\vphantom{\intercal}} f_{l}^\intercal + f_{l+n}^{\vphantom{\intercal}}f_{l+n}^\intercal\right) \\ \nonumber
& - c_{lm} \left(f_{l}^{\vphantom{\intercal}}f_{m}^\intercal + f_{m}^{\vphantom{\intercal}}f_{l}^\intercal + f_{l+n}^{\vphantom{\intercal}}f_{m+n}^\intercal + f_{m+n}^{\vphantom{\intercal}}f_{l+n}^\intercal \right) \\ \label{Zlm} & + s_{lm} \left( f_{l}^{\vphantom{\intercal}}f_{m+n}^\intercal + f_{m+n}^{\vphantom{\intercal}}f_{l}^\intercal - f_{l+n}^{\vphantom{\intercal}}f_{m}^\intercal - f_{m}^{\vphantom{\intercal}}f_{l+n}^\intercal\right) \end{align}\begin{align} %\\  
\nonumber \mathbf{Z}_{ml} & = g_{lm} \left(f_{m}^{\vphantom{\intercal}}f_{m}^\intercal + f_{m+n}^{\vphantom{\intercal}}f_{m+n}^\intercal\right) \\ \nonumber
& - c_{ml} \left(f_{l}^{\vphantom{\intercal}}f_{m}^\intercal + f_{m}^{\vphantom{\intercal}}f_{l}^\intercal + f_{l+n}^{\vphantom{\intercal}}f_{m+n}^\intercal + f_{m+n}^{\vphantom{\intercal}}f_{l+n}^\intercal \right) \\ \label{Zml} & - s_{ml} \left( f_{l+n}^{\vphantom{\intercal}}f_{m}^\intercal + f_{m}^{\vphantom{\intercal}}f_{l+n}^\intercal - f_{l}^{\vphantom{\intercal}}f_{m+n}^\intercal - f_{m+n}^{\vphantom{\intercal}}f_{l}^\intercal \right) \\ 
\nonumber \mathbf{\bar{Z}}_{ml} & = -\left(\frac{2 b_{lm} + b_{sh,lm}}{2\tau_{lm}^2} \right) \left(f_{l}^{\vphantom{\intercal}}f_{l}^\intercal + f_{l+n}^{\vphantom{\intercal}}f_{l+n}^\intercal\right) \\ \nonumber
& + c_{lm} \left( f_{l}^{\vphantom{\intercal}}f_{m+n}^\intercal + f_{m+n}^{\vphantom{\intercal}}f_{l}^\intercal - f_{l+n}^{\vphantom{\intercal}}f_{m}^\intercal - f_{m}^{\vphantom{\intercal}}f_{l+n}^\intercal\right) \\ 
 \label{Zlmbar} & + s_{lm} \left(f_{l}^{\vphantom{\intercal}}f_{m}^\intercal + f_{m}^{\vphantom{\intercal}}f_{l}^\intercal + f_{l+n}^{\vphantom{\intercal}}f_{m+n}^\intercal + f_{m+n}^{\vphantom{\intercal}}f_{l+n}^\intercal \right) \\
\nonumber \mathbf{\bar{Z}}_{ml} & = -\left(b_{lm} + \frac{b_{sh,lm}}{2}\right) \left(f_{m}^{\vphantom{\intercal}}f_{m}^\intercal + f_{m+n}^{\vphantom{\intercal}}f_{m+n}^\intercal\right) \\ \nonumber
& + c_{ml} \left( f_{l+n}^{\vphantom{\intercal}}f_{m}^\intercal + f_{m}^{\vphantom{\intercal}}f_{l+n}^\intercal - f_{l}^{\vphantom{\intercal}}f_{m+n}^\intercal - f_{m+n}^{\vphantom{\intercal}}f_{l}^\intercal \right)\\ \label{Zmlbar} & - s_{ml} \left(f_{l}^{\vphantom{\intercal}}f_{m}^\intercal + f_{m}^{\vphantom{\intercal}}f_{l}^\intercal + f_{l+n}^{\vphantom{\intercal}}f_{m+n}^\intercal + f_{m+n}^{\vphantom{\intercal}}f_{l+n}^\intercal \right).
\end{align}
\end{subequations}
Define the vector of voltage components:
\begin{equation}\label{x}
x = \begin{bmatrix} V_{d1} & V_{d2} & \ldots & V_{dn} & V_{q1} &
  V_{q2} & \ldots & V_{qn}\end{bmatrix}^\intercal
\end{equation}
and the rank-one matrix:
\begin{equation}\label{W} \mathbf{W} = x x^\intercal. \end{equation}

The active and reactive power injections at bus~$k$ are $\mathrm{tr}\left(\mathbf{Y}_k \mathbf{W} \right)$ and $\mathrm{tr}\left(\mathbf{\bar{Y}}_k \mathbf{W} \right)$, respectively, where $\mathrm{tr}\left(\cdot\right)$ indicates the matrix trace operator. The square of the voltage magnitude at bus~$k$ is $\mathrm{tr}\left(\mathbf{M}_k \mathbf{W} \right)$. The constraint $\mathrm{tr}\left(\mathbf{N}_1 \mathbf{W} \right) = 0$ sets the reference angle.

Replacing the rank-one requirement from~\eqref{W} by the less stringent constraint $\mathbf{W} \succeq 0$, where $\succeq 0$ indicates positive semidefiniteness, yields the
SDP relaxation of \eqref{opf}:

\begin{subequations}
\label{sdpprimal}
\small
\begin{align}
\label{sdp_obj} & \min_{\mathbf{W},\alpha,P_G} \sum_{k \in \mathcal{G}} \alpha_k \qquad \mathrm{subject\; to} \hspace{-20pt} & \\
\label{sdp_Pgk} & \quad P_{Gk} = \mathrm{tr}\left(\mathbf{Y}_k \mathbf{W} \right) + P_{Dk}  & \forall k \in \mathcal{N} \\
\label{sdp_P} &  \quad P_{Gk}^{\mathrm{min}} \leq P_{Gk} \leq P_{Gk}^{\mathrm{max}} & \forall k \in \mathcal{N} \\
\label{sdp_Q} &  \quad Q_{Gk}^{\mathrm{min}} \leq \mathrm{tr}\left(\mathbf{\bar{Y}}_k \mathbf{W} \right) + Q_{Dk} \leq Q_{Gk}^{\mathrm{max}} &  \forall k \in \mathcal{N} \\
\label{sdp_V} &  \quad \left(V_{k}^{\mathrm{min}}\right)^2 \leq \mathrm{tr}\left(\mathbf{M}_k \mathbf{W} \right) \leq \left(\vphantom{V_{k}^{\mathrm{min}}} V_{k}^{\mathrm{max}}\right)^2 &  \forall k \in \mathcal{N}  \\
\label{sdp_Vref} & \quad \mathrm{tr}\left(\mathbf{N}_1 \mathbf{W} \right) = 0 \\
\nonumber & \quad \left(1-c_{k1}P_{Gk}-c_{k0} + \alpha_k \right) \\ \label{sdp_quadcost} & \qquad \geq \left|\left| \begin{bmatrix} \left(1+c_{k1}P_{Gk}+c_{k0}-\alpha_k \right) \\ 2\sqrt{c_{k2}} P_{Gk} \end{bmatrix} \right|\right|_2 & \forall k \in \mathcal{G} \\
\label{sdp_Slm} & \quad S_{lm}^\mathrm{max} \geq \left|\left| \begin{bmatrix} \mathrm{tr}\left(\mathbf{Z}_{lm} \mathbf{W} \right) \\ \mathrm{tr}\left(\mathbf{\bar{Z}}_{lm} \mathbf{W} \right) \end{bmatrix} \right|\right|_2 & \forall \left(l,m \right) \in \mathcal{L} \\ %\end{align}\begin{align}
\label{sdp_Sml} & \quad S_{lm}^\mathrm{max} \geq \left|\left| \begin{bmatrix} \mathrm{tr}\left(\mathbf{Z}_{ml} \mathbf{W} \right) \\ \mathrm{tr}\left(\mathbf{\bar{Z}}_{ml} \mathbf{W} \right) \end{bmatrix} \right|\right|_2 & \forall \left(l,m \right) \in \mathcal{L} \\
\label{sdp_W} & \quad \mathbf{W} \succeq 0.
\end{align}
\end{subequations}

The generation cost constraint~\eqref{opf_obj} is implemented using the auxiliary variable $\alpha_k$ and the SOCP formulation in~\eqref{sdp_quadcost}. The apparent power line flow constraints~\eqref{opf_Slm} and \eqref{opf_Sml} are implemented with the SOCP formulations in~\eqref{sdp_Slm} and \eqref{sdp_Sml}. See~\cite{molzahn_holzer_lesieutre_demarco-large_scale_sdp_opf} for a more general formulation of the SDP relaxation that considers the possibilities of multiple generators per bus and convex piecewise-linear generation costs.

Note that rather than explicitly constraining the reference angle, the constraint~\eqref{opf_Vref} can be used to eliminate the variable $V_{q1}$ from the problem. Eliminating $V_{q1}$ removes the $n+1$ row and column from $\mathbf{W}$, with corresponding modifications to all matrices in~\eqref{sdpprimal} and removal of~\eqref{sdp_Vref}. This approach is often numerically superior to explicitly constraining the reference angle as in~\eqref{sdp_Vref}.

If the condition $\mathrm{rank}\left(\mathbf{W}\right) = 1$ is satisfied, the relaxation is ``exact'' and the global solution to~\eqref{opf} is recovered using an eigen-decomposition. Let $\lambda$ be the non-zero eigenvalue of a rank-one solution $\mathbf{W}$ to~\eqref{sdpprimal} with associated unit-length eigenvector $\eta$. The globally optimal voltage phasor is
\begin{equation}\label{Vstar}
V^\ast = \sqrt{\lambda} \left(\eta_{1:n} + \mathbf{j} \eta_{\left(n+1\right):2n} \right)
\end{equation}
where subscripts denote vector entries in MATLAB notation. 

The computational bottleneck of the SDP relaxation is the constraint~\eqref{sdp_W}, which enforces positive semidefiniteness for a $2n \times 2n$ matrix. Solving the SDP relaxation of large OPF problems requires exploiting network sparsity. A matrix completion decomposition exploits sparsity by converting the positive semidefinite constraint on the large $\mathbf{W}$ matrix~\eqref{sdp_W} to positive semidefinite constraints on many smaller submatrices of $\mathbf{W}$. These submatrices are defined using the cliques (i.e., completely connected subgraphs) of a chordal extension of the power system network graph. See~\cite{jabr11,molzahn_holzer_lesieutre_demarco-large_scale_sdp_opf,andersen2014} for a full description of a formulation that enables solution of~\eqref{sdpprimal} for systems with thousands of buses.

\section{Laplacian Objective Function}
\label{l:lap_pen}

The SDP relaxation in Section~\ref{l:sdp_relaxation} globally solves
many OPF problems~\cite{lavaei_tps,molzahn_holzer_lesieutre_demarco-large_scale_sdp_opf}. However, there are example problems for which the SDP relaxation fails to yield the globally optimal decision variables (i.e., the solution to the SDP relaxation does not satisfy the rank condition~\eqref{W}). This section proposes an approach for finding feasible points near the global optima of many problems for which the lower bounds from the SDP relaxation are close to the globally optimal objective values.

The proposed approach constrains the generation cost to be close to the lower bound obtained from the SDP relaxation. This enables the specification of an objective function based on a weighted Laplacian matrix that yields \emph{feasible} (i.e., rank-one) solutions to  many OPF problems. An iterative algorithm based on \emph{line flow mismatches} is used to determine the weights for the Laplacian matrix. %We emphasize that this approach does not require external choice of parameters (other than the maximum relaxation gap) as in related penalization techniques~\cite{lavaei_mesh,lavaei_allerton2014}.

\subsection{Generation Cost Constraint}
\label{l:gencost}

The proposed approach exploits the empirical observation that the SDP relaxation provides a very close lower bound on the optimal objective value of many typical OPF problems (i.e., there is a very small \emph{relaxation gap}). For instance, the SDP relaxation gaps for the large-scale Polish~\cite{matpower} and \mbox{PEGASE}~\cite{pegase} systems, which represent portions of European power systems, are all less than 0.3\%.\footnote{To obtain satisfactory convergence of the SDP solver, these systems are pre-processed to remove low-impedance lines (i.e., lines whose impedance values have magnitudes less than $1\times 10^{-3}$ per unit) as in~\cite{cdc2015}.}\footnote{These relaxation gaps are calculated using the objective values from the SDP relaxation~\eqref{sdpprimal} and solutions obtained either from the second-order moment relaxation~\cite{molzahn_hiskens-sparse_moment_opf} (where possible) or from \matpower{}~\cite{matpower}.} Further, the SDP relaxation is \emph{exact} (i.e., zero relaxation gap) for the IEEE \mbox{14-,} \mbox{30-,} \mbox{39-,} \mbox{57-bus} systems, the 118-bus system modified to enforce a small minimum line resistance~\cite{lavaei_tps}, and several of the large-scale Polish test cases~\cite{molzahn_holzer_lesieutre_demarco-large_scale_sdp_opf}.\footnote{Even the minor modifications performed when pre-processing low-impedance lines and enforcing minimum line resistances are not needed for some test cases. For instance, the SDP relaxation is exact for the Polish systems \mbox{2376sp,} \mbox{2737sop,} \mbox{2746wp,} and \mbox{2746wop} without modifications~\cite{molzahn_holzer_lesieutre_demarco-large_scale_sdp_opf}.} (See~\cite{matpower,ieee_test_cases} for case descriptions.) Numerical experiments also demonstrate that the SDP relaxation is exact for a variety of test cases with multiple local optima (e.g., WB2, WB3, WB5mod, and the 22- and 30-bus loop systems in~\cite{bukhsh2011}). To further demonstrate the capabilities of the SDP relaxation, 1000 modified versions were created for each of the IEEE \mbox{14-,} \mbox{30-,} \mbox{39-,} and 57-bus systems using normal random perturbations (zero-mean, 10\% standard deviation) of the load demands and power generation limits. The SDP relaxation was exact (or proved infeasibility) for 100\% and 98.7\% of the test cases derived from the \mbox{14-} and 57-bus systems, respectively. After modifications to enforce a $1\times 10^{-4}$ per unit minimum line resistance, the SDP relaxation was exact or proved infeasibility for 81.8\% and 81.2\% of the test cases derived from the \mbox{30-} and 39-bus systems, respectively.

% First-order relaxation for WB5: 946.5314
% Second-order relaxation for WB5: 946.5922

% First-order relaxation for case39mod: 2.3681e+04
% Second-order relaxation for case39mod: 2.3858e+04
% 0.74\% below the global optimum

% First-order relaxation for case3hicss:  234.3751
% Second-order relaxation for case3hicss: 234.9237

% First-order relaxation for case3iscas: 451.2952
% Second-order relaxation for case3iscas: 568.66
% 20.6\% below the global optimum

% First-order relaxation for case5hiskens: 467.0001
% Second-order relaxation for case5hiskens: 512.4540
% 8.9\% below the global optimum

% First-order relaxation for case9mod: 2.7530e+03
% Second-order relaxation for case9mod: 3.0875e+03
% 10.8\% below the global optimum

This section assumes that the lower bound provided by the SDP relaxation is within a given percentage $\delta$ of the global optimum to the OPF problem. (Most of the examples in Section~\ref{l:results} specify $\delta = 0.5\%$.) We constrain the generation cost using this assumption:
\begin{equation}\label{gencost_constraint}
\sum_{k \in \mathcal{G}} c_{k2} P_{Gk}^2 + c_{k1} P_{Gk} + c_{k0} \leq c^* \left(1 + \delta \right)
\end{equation}
where $c^*$ is the lower bound on the optimal objective value of~\eqref{opf} obtained from the semidefinite relaxation~\eqref{sdpprimal}. This constraint is implemented by augmenting the SDP relaxation's constraints \eqref{sdp_Pgk}--\eqref{sdp_W} with
\begin{equation}\label{gencost_constraint_alpha}
\sum_{k\in\mathcal{G}} \alpha_k \leq c^* \left(1 + \delta \right).
\end{equation}

If the SDP relaxation~\eqref{sdpprimal} is feasible, the feasible space defined by \eqref{sdp_Pgk}--\eqref{sdp_W} and \eqref{gencost_constraint_alpha} is non-empty for any choice of $\delta \geq 0$.\footnote{Infeasibility of the SDP relxation~\eqref{sdpprimal} certifies infeasibility of the original OPF problem~\eqref{opf}.} However, if $\delta$ is too small, there may not exist a rank-one matrix $\mathbf{W}$ (i.e., a feasible point for the original OPF problem~\eqref{opf}) in the feasible space.

The lack of a priori guarantees on the size of the relaxation gap is a challenge that the proposed approach shares with many related approaches for convex relaxations of the optimal power flow problem. Existing sufficient conditions that guarantee zero relaxation gap generally require satisfaction of non-trivial technical conditions and a limited set of network topologies~\cite{low_tutorial,lavaei_mesh}. The SDP relaxation is, however, exact for a significantly broader class of OPF problems than those that have a priori exactness guarantees, and has a small relaxation gap for an even broader class of OPF problems.\footnote{None of the aforementioned IEEE test cases, Polish systems, and \mbox{PEGASE} systems satisfy any known sufficient conditions for exactness of the SDP relaxation, but many still have zero or very small relaxation gaps.}

There are test cases that are specifically constructed to exhibit somewhat anomalous behavior in order to test the limits of the convex relaxations. The SDP relaxation gap is not small for some of these test cases. For instance, the \mbox{3-bus} system in~\cite{iscas2015}, the \mbox{5-bus} system in~\cite{hicss2014}, and the \mbox{9-bus} system in~\cite{bukhsh_tps} have relaxation gaps of 20.6\%, 8.9\%, and 10.8\%, respectively, and the test cases in~\cite{kocuk15} have relaxation gaps as large as 52.7\%. The approach proposed in this paper is not appropriate for such problems. Future progress in convex relaxation theory is required to develop broader conditions that provide a priori certification that the SDP relaxation is exact or has a small relaxation gap. We also await the development of more extensive sets of OPF test cases to further explore the observation that many typical existing practical test cases have small SDP relaxation gaps.

\subsection{Laplacian Objective Function}

Consider the optimization problem
\begin{equation}
\label{augmented_sdp}
\begin{aligned}
& \min_{\mathbf{W},\alpha,P_G} f\left(\mathbf{W}\right) \\ & \mathrm{subject\; to\;\;} \eqref{sdp_Pgk}-\eqref{sdp_W},\; \eqref{gencost_constraint_alpha}
\end{aligned}
\end{equation}
where $f\left(\mathbf{W} \right)$ is an arbitrary linear function. Any solution to~\eqref{augmented_sdp} with $\mathrm{rank}\left(\mathbf{W}\right) = 1$ yields a feasible solution to the OPF problem~\eqref{opf} within $\delta$ of the globally optimal objective value due to the constraint \eqref{gencost_constraint_alpha} on the generation cost. This constraint effectively frees the choice of the function $f\left(\mathbf{W} \right)$ to obtain a \emph{feasible} rather than \emph{minimum-cost} solution to \eqref{opf}.

We therefore seek an objective function $f\left(\mathbf{W}\right)$
which maximizes the likelihood of obtaining
$\mathrm{rank}\left( \mathbf{W}\right) = 1$. This section describes a
Laplacian form for the function
$f\left(\mathbf{W}\right)$. Specifically, we consider a $n_l \times
n_l$ diagonal matrix $\mathbf{D}$ containing weights for the network
Laplacian matrix $\mathbf{L} = \mathbf{A}_{inc}^\intercal \mathbf{D}
\mathbf{A}_{inc}$, where $\mathbf{A}_{inc}$ is the $n_l\times n$
incidence matrix for the network. The off-diagonal term
$\mathbf{L}_{ij}$ is equal to the negative of the sum of the weights
for the lines connecting buses~$i$ and~$j$, and the diagonal term
$\mathbf{L}_{ii}$ is equal to the sum of the weights of the lines
connected to bus~$i$. The objective function is
\begin{equation}\label{lap_obj}
f\left(\mathbf{W}\right) = \mathrm{tr}\left(\begin{bmatrix}\mathbf{L} & \mathbf{0}_{n\times n} \\ \mathbf{0}_{n\times n} & \mathbf{L}\end{bmatrix}\mathbf{W} \right).
\end{equation}

The choice of an objective function based on a Laplacian matrix is
motivated by previous literature. An existing penalization
approach~\cite{lavaei_mesh} augments the objective
  function by adding a term that minimizes the total reactive power
injection. This reactive power penalty can
be implemented by adding the term
\begin{equation}\label{Qpen}
\epsilon_b\, \mathrm{tr}\left(
\begin{bmatrix}\hphantom{-}\Re\left(\frac{\mathbf{Y}^H - \mathbf{Y}}{2\mathbf{j}}\right) & \Im\left(\frac{\mathbf{Y}^H - \mathbf{Y}}{2\mathbf{j}}\right) \\
-\Im\left(\frac{\mathbf{Y}^H - \mathbf{Y}}{2\mathbf{j}}\right) & \Re\left(\frac{\mathbf{Y}^H - \mathbf{Y}}{2\mathbf{j}}\right)
\end{bmatrix} \mathbf{W}\right)
\end{equation}
to the objective function of the SDP relaxation~\eqref{sdp_obj}, where
$\epsilon_b$ is a specified penalty parameter and
$\left(\cdot\right)^H$ indicates the complex conjugate transpose
operator. In the absence of phase-shifting transformers (i.e.,
$\theta_{lm} = 0 \quad \forall \left(l,m\right) \in \mathcal{L}$), the
matrix $ \frac{\mathbf{Y}^H - \mathbf{Y}}{2\mathbf{j}}$ is equivalent
to $-\Im\left(\mathbf{Y}\right) = -\mathbf{B}$, which is a
\emph{weighted Laplacian matrix} (with weights determined by the
branch susceptance parameters $b_{lm} = \frac{-X_{lm}}{R_{lm}^2 +
  X_{lm}^2}$) plus a diagonal matrix composed of shunt susceptances.

Early work on SDP relaxations of OPF problems~\cite{lavaei_tps} advocates enforcing a minimum resistance of $\epsilon_r$ for all lines in the network. For instance, the SDP relaxation fails to be exact for the IEEE 118-bus system~\cite{ieee_test_cases}, but the relaxation is exact after enforcing a minimum line resistance of $\epsilon_r = 1\times 10^{-4}$~per~unit. After enforcing a minimum line resistance, the active power losses are given by
\begin{equation}
\mathrm{tr}\left(\begin{bmatrix}
\hphantom{-}\Re\left(\frac{\mathbf{Y}_r + \mathbf{Y}_r^H}{2}\right) & \Im\left(\frac{\mathbf{Y}_r + \mathbf{Y}_r^H}{2}\right) \\
-\Im\left(\frac{\mathbf{Y}_r + \mathbf{Y}_r^H}{2}\right) & \Re\left(\frac{\mathbf{Y}_r + \mathbf{Y}_r^H}{2}\right) \end{bmatrix}\mathbf{W}\right)
\end{equation}
where $\mathbf{Y}_r$ is the network admittance matrix after enforcing a minimum branch resistance of $\epsilon_r$. In the absence of phase-shifting transformers, $\frac{\mathbf{Y}_r + \mathbf{Y}_r^H}{2}$ is equivalent to $\Re\left(\mathbf{Y}_r\right)$, which is a \emph{weighted Laplacian matrix} (with weights determined by the branch conductance parameters $g_{lm} = \frac{R_{lm}}{R_{lm}^2 + X_{lm}^2}$) plus a diagonal matrix composed of shunt conductances. Since typical OPF problems have objective functions that increase with active power losses, enforcing minimum line resistances is similar to a weighted Laplacian penalization.\footnote{Note that since enforcing minimum line resistances also affects the power injections and the line flows, the minimum line resistance cannot be solely represented as a Laplacian penalization of the objective function.}

% \mathbf{A}_{inc}^\intercal\mathbf{D}_r\mathbf{A}_{inc}
% $\mathbf{D}_r$ is an $n_l\times n_l$ diagonal matrix with diagonal elements corresponding to lines with $R_{lm} < \epsilon_r$ having values of $\epsilon_r - R_{lm}$ with the remaining elements equal to zero.

% where $\Delta \mathbf{Y} = \mathbf{Y}_1 - \mathbf{Y}_0$ is the difference between the admittance matrix of the original network, which is denoted $\mathbf{Y}_0$, and the admittance matrix after enforcing a minimum line resistance of $\epsilon_r$, which is denoted $\mathbf{Y}_1$.

The proposed objective function~\eqref{lap_obj} is equivalent to a linear combination of certain components of $\mathbf{W}$:
\begin{align}\nonumber
f\left(\mathbf{W}\right) = \sum_{\left(l,m\right) \in \mathcal{L}}&  \mathbf{D}_{\left(l,m\right)} \left(\mathbf{W}_{ll} - 2\mathbf{W}_{lm} + \mathbf{W}_{mm} \right. \\
&  \left. + \mathbf{W}_{l+n,l+n} - 2\mathbf{W}_{l+n,m+n} + \mathbf{W}_{m+n,m+n} \right)
\end{align}
where $\mathbf{D}_{\left(l,m\right)}$ is the diagonal element of $\mathbf{D}$ corresponding to the line from bus~$l$ to bus~$m$. If $\mathbf{W} = xx^\intercal$ (i.e., $\mathbf{W}$ is a rank-one matrix) with $x$ defined as in~\eqref{x}, then the objective function \eqref{lap_obj} is equivalent to 
\begin{align}\nonumber
f\left(\mathbf{xx^\intercal}\right) & = \sum_{\left(l,m\right) \in \mathcal{L}} \mathbf{D}_{\left(l,m\right)} \left\lbrace \left(x_l - x_m \right)^2 + \left(x_{l+n} - x_{m+n} \right)^2 \right\rbrace \\ \label{lap_pen_xxt}
& = \sum_{\left(l,m\right) \in \mathcal{L}}  \mathbf{D}_{\left(l,m\right)} \left\lbrace \left(V_{dl} - V_{dm} \right)^2 + \left(V_{ql} - V_{qm} \right)^2 \right\rbrace.
\end{align}
Note that the Laplacian objective function~\eqref{lap_pen_xxt} is convex in the voltage components $V_d$ and $V_q$ when the weights in $\mathbf{D}$ are non-negative. This is in contrast to the reactive power penalization in~\cite{lavaei_mesh}: an objective function that penalizes reactive power injections may be non-convex in terms of the voltage components $V_d$ and $V_q$ when the network has non-zero shunt capacitors.

When the SDP relaxation fails to yield the global optimum,
  the relaxation often ``artificially'' increases the voltage
  magnitudes to reduce active power losses. This results in voltage
  magnitudes and power injections that are feasible for the
  relaxation~\eqref{sdpprimal} but infeasible for the OPF
  problem~\eqref{opf}. By minimizing the squared differences between
  the voltage phasors at connected buses, the Laplacian objective
  function counteracts this tendency of the SDP
  relaxation. Intuitively, the proposed approach uses the Laplacian
  objective function to balance two potentially competing tendencies:
  increasing voltage magnitudes to reduce active power losses such
  that the generation cost constraint is satisfied
  versus decreasing voltage differences to reduce the
  Laplacian objective function.

From a physical perspective, the Laplacian objective's tendency to
reduce voltage differences is similar to both the reactive power
penalization proposed in~\cite{lavaei_mesh} and the minimum branch
resistance advocated in~\cite{lavaei_tps}. For typical operating
conditions, reactive power injections are closely related to voltage
magnitude differences, so penalizing reactive power injections tends
to result in solutions with similar voltages. Likewise,
the active power losses associated with line resistances increase with
the square of the current flow through the line, which is determined
by the voltage difference across the line. Thus, enforcing a minimum
line resistance tends to result in solutions with smaller voltage
differences in order to reduce losses.

In addition, a Laplacian regularizing term has been used to obtain desirable solution characteristics for a variety of other optimization problems (e.g., machine learning problems~\cite{smola2003,melacci2011}, sensor network localization problems~\cite{weinberger2006}, and analyses of flow networks~\cite{taylor2011}).

%With non-negative weights in $\mathbf{D}$ and in the absence of any constraints, the objective function \eqref{lap_obj} achieves a minimum value of zero with all voltage phasors taking the same value.

\subsection{An Algorithm for Determining the Laplacian Weights}

Having established a weighted Laplacian form for the objective function, we introduce an iterative algorithm for determining appropriate weights $\mathbf{D}$ for obtaining a solution to~\eqref{augmented_sdp} with $\mathrm{rank}\left(\mathbf{W} \right) = 1$. We note that the proposed algorithm is similar in spirit to the method in~\cite[Section~2.4]{candes-2013}, which iteratively updates weighting parameters to promote low-rank solutions of SDPs related to image reconstruction problems.

The proposed algorithm is inspired by the apparent power line flow
penalty used in~\cite{lavaei_allerton2014} and the iterative approach
to determining appropriate buses for enforcing higher-order moment
constraints in~\cite{molzahn_hiskens-sparse_moment_opf}. The approach
in~\cite{lavaei_allerton2014} penalizes the apparent power flows on
lines associated with certain submatrices of $\mathbf{W}$ that are not
rank one.\footnote{The submatrices are determined by the maximal
  cliques of a chordal supergraph of the network;
  see~\cite{jabr11,molzahn_holzer_lesieutre_demarco-large_scale_sdp_opf,lavaei_allerton2014}
  for further details.} Similar to the approach
in~\cite{lavaei_allerton2014}, the proposed algorithm adds terms to
the objective function that are associated with certain
``problematic lines.''

The heuristic for identifying problematic lines is inspired by the approach used in~\cite{molzahn_hiskens-sparse_moment_opf} to detect ``problematic buses'' for application of higher-order moment constraints. Denote the solution to~\eqref{augmented_sdp} as $\mathbf{W}^\star$ and the closest rank-one matrix to $\mathbf{W}^\star$ as $\mathbf{W}^{\left(1\right)}$. (By the Eckart and Young theorem~\cite{eckart36}, the closest rank-one matrix is calculated using the eigendecomposition $\mathbf{W}^{\left(1\right)} = \lambda_1 \eta_1 \eta_1^\intercal$, where $\lambda_1$ and $\eta_1$ are the largest eigenvalue of $\mathbf{W}^\star$ and its associated unit-length eigenvector, respectively.) If $\mathbf{W}^\star = \mathbf{W}^{\left(1\right)}$, then $\mathrm{rank}\left(\mathbf{W}^\star\right) = 1$ and we can recover the global optimum to~\eqref{opf} using~\eqref{Vstar}. Otherwise, previous work~\cite{molzahn_hiskens-sparse_moment_opf} compares the power injections associated with $\mathbf{W}^\star$ and $\mathbf{W}^{\left(1\right)}$ to calculate power injection mismatches $S_k^{inj\,mis}$ for each bus~$k \in \mathcal{N}$:
\begin{align}\nonumber
& S_{k}^{inj\,mis} = \\ \label{sdpmis} & \;\; \left| 
\mathrm{tr}\left\lbrace\mathbf{Y}_k \left(\mathbf{W}^\star-\mathbf{W}^{\left(1\right)}\right)\right\rbrace  + \mathbf{j}
\mathrm{tr}\left\lbrace\mathbf{\bar{Y}}_k \left(\mathbf{W}^\star-\mathbf{W}^{\left(1\right)}\right)\right\rbrace\right|
\end{align}
where $\left|\,\cdot\, \right|$ denotes the magnitude of the complex argument. In the parlance of~\cite{molzahn_hiskens-sparse_moment_opf}, problematic buses are those with large power injection mismatches $S_{k}^{inj\,mis}$.

To identify problematic lines rather than buses, we modify~\eqref{sdpmis} to calculate apparent power flow mismatches $S_{\left(l,m\right)}^{flow\,mis}$ for each line~$\left(l,m\right) \in \mathcal{L}$:
\begin{align}\nonumber
& S_{\left(l,m\right)}^{flow\,mis} = \\ \nonumber &\;\; \left| 
\mathrm{tr}\left\lbrace\mathbf{Z}_{lm} \left(\mathbf{W}^\star-\mathbf{W}^{\left(1\right)}\right)\right\rbrace + \mathbf{j}
\mathrm{tr}\left\lbrace\mathbf{\bar{Z}}_{lm} \left(\mathbf{W}^\star-\mathbf{W}^{\left(1\right)}\right)\right\rbrace\right| \\ \label{flowmis}
&\;\; + \left| 
\mathrm{tr}\left\lbrace\mathbf{Z}_{ml} \left(\mathbf{W}^\star-\mathbf{W}^{\left(1\right)}\right)\right\rbrace + \mathbf{j}
\mathrm{tr}\left\lbrace\mathbf{\bar{Z}}_{ml} \left(\mathbf{W}^\star-\mathbf{W}^{\left(1\right)}\right)\right\rbrace\right|.
\end{align}
Observe that $S_{\left(l,m\right)}^{flow\,mis}$ sums the magnitude of the apparent power flow mismatches at both ends of each line.

The condition $\mathrm{rank}\left(\mathbf{W}^\star\right) = 1$ (i.e., ``feasibility'' in this context) is considered satisfied for practical purposes using the criterion that the maximum line flow and power injection mismatches (i.e., $\max_{\left(l,m\right)\in\mathcal{L}} S_{\left(l,m\right)}^{flow\,mis}$ and $\max_{k\in\mathcal{N}}S_{k}^{inj\,mis}$) are less than specified tolerances~$\epsilon_{flow}$ and~$\epsilon_{inj}$, respectively, and the voltage magnitude limits~\eqref{opf_V} are satisfied to within a specified tolerance~$\epsilon_{V}$.\footnote{For all test cases, the voltage magnitude limits were satisfied whenever the power injection and line flow mismatch tolerances were achieved.}

\newfloat{algorithm}{bt}{lop}
\begin{algorithm}
\caption{Iterative Algorithm for Determining Weights}\label{a:weights}
\begin{algorithmic}[1]
\State \textbf{Input}: tolerances $\epsilon_{flow}$ and $\epsilon_{inj}$, max relaxation gap $\delta$
\State Set $\mathbf{D} = \mathbf{0}_{n_l\times n_l}$
\State Solve the SDP relaxation~\eqref{sdpprimal} to obtain $c^*$
\State Calculate $S^{flow\,mis}$ and $S^{inj\,mis}$ using~\eqref{flowmis} and~\eqref{sdpmis}
\NoDo
\While{termination criteria not satisfied} % {\small $\max\left\lbrace S^{flow\,mis} \right\rbrace \geq \epsilon_{flow}$ or $\max\left\lbrace S^{inj\,mis} \right\rbrace \geq \epsilon_{inj}$}
	\State Update weights: $\mathbf{D} \leftarrow \mathbf{D} + \mathrm{diag}\left(S^{flow\,mis}\right)$
	\State Solve the generation-cost-constrained relaxation~\eqref{augmented_sdp}
	\State Calculate $S^{flow\,mis}$ and $S^{inj\,mis}$ using~\eqref{flowmis} and~\eqref{sdpmis}
\EndWhile
\State Calculate the voltage phasors using~\eqref{Vstar} and terminate
\end{algorithmic}
\end{algorithm}

As described in Algorithm~\ref{a:weights}, the weights on the diagonal of $\mathbf{D}$ are determined from the line flow mismatches $S_{\left(l,m\right)}^{flow\,mis}$. Specifically, the proposed algorithm first solves the SDP relaxation~\eqref{sdpprimal} to obtain both the lower bound $c^*$ on the optimal objective value and the initial line flow and power injection mismatches $S_{\left(l,m\right)}^{flow\,mis},\, \forall \left(l,m\right) \in\mathcal{L}$ and $S_{k}^{inj\,mis},\, \forall k\in\mathcal{N}$. 

While the termination criteria ($\max_{\left(l,m\right)\in\mathcal{L}} \big\{ S_{\left(l,m\right)}^{flow\,mis} \big\} <
\epsilon_{flow}$, $\max_{k\in\mathcal{N}}  \big\{
  S_{k}^{inj\,mis} \big\} < \epsilon_{inj}$, and no
  voltage limits violated by more than $\epsilon_V$) are not
satisfied, the algorithm solves~\eqref{augmented_sdp} (i.e., the SDP
relaxation with the constraint ensuring that the generation cost is
within $\delta$ of the lower bound). The objective function is defined
using the weighting matrix $\mathbf{D} =
\mathrm{diag}\left(S^{flow\,mis}\right)$, where
$\mathrm{diag}\left(\cdot\right)$ denotes the matrix with the vector
argument on the diagonal and other entries equal to zero. Each
iteration adds the line flow mismatch vector $S^{flow\,mis}$ from the
solution to~\eqref{augmented_sdp} to the previous weights (i.e.,
$\mathbf{D} \leftarrow \mathbf{D} +
\mathrm{diag}\left(S^{flow\,mis}\right)$).

Upon satisfaction of the termination criteria, the algorithm uses~\eqref{Vstar} to recover a feasible solution to~\eqref{opf} that has an objective value within $\delta$ of the global optimum. Again, ``feasibility'' in this context is judged using the termination criteria $\epsilon_{flow}$, $\epsilon_{inj}$, and $\epsilon_V$.

Note that Algorithm~\ref{a:weights} is not guaranteed to converge. Non-convergence may be due to the value of $\delta$ being too small (i.e., there does not exist a rank-one solution that satisfies~\eqref{gencost_constraint_alpha}) or failure to find a rank-one solution that does exist. To address the former case, Algorithm~\ref{a:weights} could be modified to include an ``outer loop'' that increments $\delta$ by a specified amount (e.g., $0.5$\%) if convergence is not achieved in a certain  number of iterations. We note that, like other convex relaxation methods, the proposed approach would benefit from further theoretical work regarding the development of a priori guarantees on the size of the relaxation gap for various classes of OPF problems.
% Increasing $\delta$ was generally unnecessary for the test cases considered in this paper.}\footnote{The exception was the test case MH118L~\cite{molzahn_hiskens-sparse_moment_opf}, which had a relaxation gap of $0.76$\% and would therefore benefit from increasing $\delta$ when no solution was obtained for $\delta = 0.5\%$. We simply used $\delta = 1\%$ for this test case.}

For some problems with large relaxation gaps (e.g., the
  \mbox{3-bus} system in~\cite{iscas2015}, the 5-bus system
  in~\cite{hicss2014}, and the 9-bus system in~\cite{bukhsh_tps}), no
  purely penalization-based methods have so far successfully addressed
  the latter case where the proposed algorithm fails to
  find a rank-one solution that satisfies the generation cost
  constraint~\eqref{gencost_constraint_alpha} with 
    sufficiently large $\delta$ (i.e., no known penalty parameters
  yield feasible solutions using the methods
  in~\cite{lavaei_mesh,lavaei_allerton2014} for these test cases). One
  possible approach for addressing this latter case is the combination
  of penalization techniques with Lasserre's moment relaxation
  hierarchy~\cite{pscc2014,cedric_msdp,ibm_paper,molzahn_hiskens-sparse_moment_opf}. The
  combination of the moment relaxations with the penalization methods
  enables the computation of near-globally-optimal solutions for a
  broader class of OPF problems than either method achieves
  individually. See~\cite{cdc2015} for further details on this
  approach.

We note that despite the lack of a convergence guarantee, the examples in Section~\ref{l:results} demonstrate that Algorithm~\ref{a:weights} is capable of finding feasible points that are near the global optimum for many OPF problems, including large test cases.

% Can we say that the closest rank one matrix gives a lower bound on the relaxation gap? Probably easiest to not exploit sparsity when doing this calculation so that the closest rank one matrix is unambigious. <-- I don't think this is actually correct. The closest rank one matrix is closest in the sense of a Frobenius norm, but not necessarily closest in terms of the cost function. In other words, it may be the case that another rank one matrix which is farther away in terms of the Frobenius norm actually has a lower cost.

\section{Results}
\label{l:results}

This section demonstrates the effectiveness of the
proposed approach using several small example problems as well as
large test cases representing portions of European power
systems. The SDP relaxation yields a small but non-zero
  relaxation gap for the test cases selected in this section, and
  Algorithm~\ref{a:weights} yields points that are feasible
  for~\eqref{opf} (to within the specified termination
    criteria) and that are near the global optimum for these test
  cases. For other test cases with a large SDP relaxation gap, such as those mentioned earlier in~\cite{bukhsh_tps,hicss2014,iscas2015,kocuk15}, the proposed algorithm does
  not converge when tested with a variety of values for $\delta$.

The results in this section use line flow and power injection mismatch tolerances $\epsilon_{flow}$ and $\epsilon_{inj}$ that are both equal to $1$~MVA and $\epsilon_V = 5\times 10^{-4}$~per~unit. The implementation of Algorithm~\ref{a:weights} uses \mbox{MATLAB 2013a}, YALMIP \mbox{2015.02.04}~\cite{yalmip}, and Mosek \mbox{7.1.0.28}~\cite{mosek}, and was solved using a computer with a quad-core 2.70~GHz processor and 16~GB of RAM.

Applying Algorithm~\ref{a:weights} to several small- to medium-size
test cases
from~\cite{allerton2011,hicss2014,molzahn_hiskens-sparse_moment_opf,ieee_test_cases,kocuk15}
yields the results shown in
Table~\ref{t:small_results}. Tables~\ref{t:large_gencost_results}
and~\ref{t:large_loss_results} show the results from applying
Algorithm~\ref{a:weights} to large test cases which minimize
generation cost and active power losses, respectively. These test
cases, which are from~\cite{matpower} and~\cite{pegase}, represent
portions of European power systems. The SDP
relaxation~\eqref{sdpprimal} has a small but non-zero relaxation gap
for all test cases considered in this section. The columns of
Tables~\ref{t:small_results}--\ref{t:large_loss_results} show the case
name and reference, the number of iterations of
Algorithm~\ref{a:weights}, the final maximum apparent power flow
mismatch $\max_{\left(l,m\right)\in\mathcal{L}} \big\{
  S_{\left(l,m\right)}^{flow\,mis} \big\}$ in MVA, the final maximum
power injection mismatch $\max_{k\in\mathcal{N}}  \big\{
  S_{k}^{inj\,mis} \big\}$ in MVA, the specified value of $\delta$, an
  upper bound on the relaxation gap from the solution to the SDP
  relaxation~\eqref{sdpprimal}, and the total solver time in seconds.

Note that the large test cases in Tables~\ref{t:large_gencost_results}
and~\ref{t:large_loss_results} were preprocessed to remove
low-impedance lines as described in~\cite{cdc2015} in order to improve
the numerical convergence of the SDP relaxation. Lines
  which have impedance magnitudes less than a threshold
  (\texttt{thrshz} in~\cite{cdc2015}) of $1\times 10^{-3}$~per~unit
  are eliminated by merging the terminal buses. Table~\ref{tab:size}
  describes the number of buses and lines before and after this
  preprocessing. Low-impedance line preprocessing was not needed for
  the test cases in Table~\ref{t:small_results}. After preprocessing,
  MOSEK's SDP solver converged with sufficient accuracy to yield
  solutions that satisfied the voltage magnitude limits to within
  \mbox{$\epsilon_V = 1\times 10^{-4}$}~per~unit and the 
    power injection and line flow constraints to within the
  corresponding mismatches shown in
  Tables~\ref{t:small_results}--\ref{t:large_loss_results}.

%\multicolumn{1}{|c|}{\textbf{Case}}  & \multicolumn{1}{c|}{\textbf{Num.}} & \multicolumn{1}{c|}{\textbf{Num.}} & \multicolumn{2}{c|}{M{\sc atpower} \textbf{Solution~\cite{matpower}}}\\ \cline{4-5} 
%\multicolumn{1}{|c|}{\textbf{Name}}  & \multicolumn{1}{c|}{\textbf{Buses}} & \multicolumn{1}{c|}{\textbf{Lines}} & \multicolumn{1}{c|}{\textbf{Gen. Cost}} & \multicolumn{1}{c|}{\textbf{Loss Min.}}\\

%Also note that the objective function should be appropriately normalized relative to the constraint coefficients in order to avoid numeric inaccuracies with the SDP solver.

% We next discuss several observations from the results in Tables~\ref{t:small_results}--\ref{t:large_loss_results}. 

These results show that Algorithm~\ref{a:weights} finds feasible points (within the specified tolerances) that have objective values near the global optimum for a variety of test cases.  Further, Algorithm~\ref{a:weights} globally solves all OPF problems for which the SDP relaxation~\eqref{sdpprimal} is exact (e.g., many of the IEEE test cases~\cite{lavaei_tps}, several of the Polish test systems~\cite{molzahn_holzer_lesieutre_demarco-large_scale_sdp_opf}, and the 89-bus PEGASE system~\cite{pegase}). Thus, the algorithm is a practical approach for addressing a broad class of OPF problems. 

We note, however, that Algorithm~\ref{a:weights} does not yield a feasible point for all OPF problems. For instance, the test case WB39mod from~\cite{bukhsh_tps} has line flow and power injection mismatches of 18.22~MVA and 12.99~MVA, respectively, after 1000 iterations of Algorithm~\ref{a:weights}. The challenge associated with this case seems to result from light loading with limited ability to absorb a surplus of reactive power injections, yielding at least two local solutions. In addition to challenging the method proposed in this paper, no known penalty parameters yield feasible solutions to this problem. Generalizations of the SDP relaxation using the Lasserre hierarchy have successfully calculated the global solution to this case~\cite{molzahn_hiskens-sparse_moment_opf,cdc2015}. Further, while Algorithm~\ref{a:weights} converges for five of the seven test cases in~\cite{kocuk15} which have small relaxation gaps (less than 2.5\%), the algorithm fails for two other such test cases as well as several other test cases in~\cite{kocuk15} which have large relaxation gaps. We note that the tree topologies used in the test cases in~\cite{kocuk15} are a significant departure from the mesh networks used in the standard test cases from which they were derived; the proposed algorithm succeeds for several test cases that share the original network topologies.

% case39mod
% For delta = 0.05, 18.22 MVA flowmis, 12.99 MVA power injection mismatch
% For delta = 0.10, 20.94 MVA flowmis, 12.86 MVA power injection mismatch
% matpower's MIPS finds (global) solution: 23858.16
% first-order real sdp: 2.370431162814909e+04
% second-order msdp: 2.385803605859051e+04

\begin{table}[t]
\centering
\caption{Results for Small and Medium Size Test Cases}
\label{t:small_results}
\resizebox{\columnwidth}{!}{%
\begin{tabular}{|l|c|c|c|c|c|c|}
\hline 
\multicolumn{1}{|c|}{\textbf{Case}} & \!\!\textbf{Num.}\!\! & \textbf{Max} & \textbf{Max} & \textbf{$\delta$} & \textbf{Max} & \!\!\textbf{Solver}\!\!\\
\multicolumn{1}{|c|}{\textbf{Name}} & \!\!\textbf{Iter.}\!\! & \!\!\textbf{Flow Mis.}\!\! & \!\!\textbf{Inj. Mis.}\!\! & \textbf{(\%)} & \!\!\textbf{Relax.}\!\! & \textbf{Time} \\
& & \textbf{(MVA)} & \textbf{(MVA)} & & \!\!\textbf{Gap (\%)}\!\! & \textbf{(sec)} \\ \hline\hline
\!\!\! LMBD3~\cite{allerton2011}\!\!\! & 1 & $1.3${\textsc e}${-5}$ & $1.6${\textsc e}${-5}$ & 0.5 & 0.50 & 0.7 \\ \hline
\!\!\! MLD3~\cite{hicss2014}\!\!\! & 1 & $7.3${\textsc e}${-6}$ & $7.2${\textsc e}${-5}$ & 0.5 & 0.50 & 0.5 \\ \hline
\!\!\! MH14Q~\cite{molzahn_hiskens-sparse_moment_opf}\!\!\! & 2 & $1.8${\textsc e}${-5}$ & $9.9${\textsc e}${-6}$ & 0.5 & 0.02 & 1.2 \\ \hline
\!\!\! MH14L~\cite{molzahn_hiskens-sparse_moment_opf}\!\!\! & 2 & $8.1${\textsc e}${-5}$ & $7.8${\textsc e}${-5}$ & 0.5 & 0.33 & 1.2 \\ \hline
\!\!\! KDS14Lin~\cite{kocuk15}\!\!\! & 1 & $1.2${\textsc e}${-3}$ & $9.2${\textsc e}${-4}$ & 1.0 & 1.00 & 0.7 \\ \hline
\!\!\! KDS14Quad~\cite{kocuk15}\!\!\! & 1 & $1.4${\textsc e}${-4}$ & $8.4${\textsc e}${-5}$ & 1.0 & 1.00 & 0.6 \\ \hline
\!\!\! KDS30Lin~\cite{kocuk15}\!\!\! & 7 & $9.3${\textsc e}${-1}$ & $9.2${\textsc e}${-1}$ & 2.5 & 2.50 & 4.6 \\ \hline
\!\!\! KDS30Quad~\cite{kocuk15}\!\!\! & 6 & $8.1${\textsc e}${-1}$ & $8.0${\textsc e}${-1}$ & 2.0 & 2.00 & 3.6 \\ \hline
\!\!\! KDS30IEEEQuad~\cite{kocuk15}\!\!\! & 100 & $9.5${\textsc e}${-1}$ & $7.2${\textsc e}${-1}$ & 2.5 & 2.50 & 129.8 \\ \hline
\!\!\! MH39L~\cite{molzahn_hiskens-sparse_moment_opf}\!\!\! & 1 & $1.3${\textsc e}${-2}$ & $9.8${\textsc e}${-3}$ & 0.5 & 0.27 & 0.7 \\ \hline
\!\!\! MH57Q~\cite{molzahn_hiskens-sparse_moment_opf}\!\!\! & 1 & $1.2${\textsc e}${-3}$ & $6.9${\textsc e}${-4}$ & 0.5 & 0.03 & 0.7 \\ \hline
\!\!\! MH57L~\cite{molzahn_hiskens-sparse_moment_opf}\!\!\! & 1 & $3.2${\textsc e}${-4}$ & $5.2${\textsc e}${-4}$ & 0.5 & 0.16 & 0.9 \\ \hline
\!\!\! MH118Q~\cite{molzahn_hiskens-sparse_moment_opf}\!\!\! & 2 & $3.3${\textsc e}${-3}$ & $2.7${\textsc e}${-3}$ & 0.5 & 0.50 & 2.6 \\ \hline
\!\!\! MH118L~\cite{molzahn_hiskens-sparse_moment_opf}\!\!\! & 2 & $3.1${\textsc e}${-3}$ & $3.1${\textsc e}${-3}$ & 1.0 & 1.00 & 3.3\\ \hline 
\!\!\! IEEE 300~\cite{ieee_test_cases}\!\!\! & 1 & $1.3${\textsc e}${-1}$ & $1.2${\textsc e}${-1}$ & 0.5 & 0.01 & 3.0\\ \hline
\end{tabular}
}
\end{table}
% \tablefootnote{The test case MH118L uses a maximum relaxation gap of $\delta = 1\%$. All other test cases use $\delta = 0.5\%$.}

%\!\!\! KDS57L~\cite{kocuk15}\!\!\! & 9 & $8.5${\textsc e}${-1}$ & $8.2${\textsc e}${-1}$ & 17.5 & 17.50 & 4.1 \\ \hline
%\!\!\! KDS57Q~\cite{kocuk15}\!\!\! & 2 & $4.5${\textsc e}${-1}$ & $3.4${\textsc e}${-1}$ & 14.0 & 17.50 & 4.1 \\ \hline

\begin{table}[t]
\centering
\caption{Results for Large Test Cases that Minimize Generation Cost}
\label{t:large_gencost_results}
\resizebox{\columnwidth}{!}{%
\begin{tabular}{|l|c|c|c|c|c|c|}
\hline 
\multicolumn{1}{|c|}{\textbf{Case}} & \!\!\textbf{Num.}\!\! & \textbf{Max} & \textbf{Max} & \textbf{$\delta$} & \textbf{Max} & \!\!\textbf{Solver}\!\!\\
\multicolumn{1}{|c|}{\textbf{Name}} & \!\!\textbf{Iter.}\!\! & \!\!\textbf{Flow Mis.}\!\! & \!\!\textbf{Inj. Mis.}\!\! & \textbf{(\%)} & \!\!\textbf{Relax.}\!\! & \textbf{Time} \\
& & \textbf{(MVA)} & \textbf{(MVA)} & & \!\!\textbf{Gap (\%)}\!\! & \textbf{(sec)} \\ \hline\hline
\!\!\! PL-2383wp~\cite{matpower}\!\!\! & 2 & 0.54 & 0.50 & 0.5 & 0.50 & 78.6 \\ \hline
\!\!\! PL-3012wp~\cite{matpower}\!\!\! & 2 & 0.36 & 0.27 & 0.5 & 0.50 & 107.6 \\ \hline
\!\!\! PL-3120sp~\cite{matpower}\!\!\! & 2 & 0.56 & 0.33 & 0.5 & 0.50 & 84.2 \\ \hline
\end{tabular}
}
\end{table}

\begin{table}[t!]
\centering
\caption{Results for Large Test Cases that Minimize Active Power Loss}
\label{t:large_loss_results}
\resizebox{\columnwidth}{!}{%
\begin{tabular}{|l|c|c|c|c|c|c|}
\hline 
\multicolumn{1}{|c|}{\textbf{Case}} & \!\!\textbf{Num.}\!\! & \textbf{Max} & \textbf{Max} & \textbf{$\delta$} & \textbf{Max} & \!\!\textbf{Solver}\!\!\\
\multicolumn{1}{|c|}{\textbf{Name}} & \!\!\textbf{Iter.}\!\! & \!\!\textbf{Flow Mis.}\!\! & \!\!\textbf{Inj. Mis.}\!\! & \textbf{(\%)} & \!\!\textbf{Relax.}\!\! & \textbf{Time} \\
& & \textbf{(MVA)} & \textbf{(MVA)} & & \!\!\textbf{Gap (\%)}\!\! & \textbf{(sec)} \\ \hline\hline
\!\!\! PL-2383wp~\cite{matpower} & 5 & 0.21 & 0.16 & 0.5 & 0.26 & 154.0 \\ \hline
\!\!\! PL-3012wp~\cite{matpower} & 5 & 0.08 & 0.04 & 0.5 & 0.18 & 232.2 \\ \hline
\!\!\! PL-3120sp~\cite{matpower} & 5 & 0.25 & 0.19 & 0.5 & 0.38 & 232.6 \\ \hline
\!\!\! PEGASE-1354~\cite{pegase}\!\!\! & 12 & 0.27 & 0.18 & 0.5 & 0.15 & 199.2 \\ \hline
\!\!\! PEGASE-2869~\cite{pegase}\!\!\! & 38 & 0.91 & 0.69 & 0.5 & 0.15 & 2378.4 \\ \hline
\end{tabular}
}
\end{table}

\begin{table}[t!]
\centering
\caption{Descriptions of Large Test Cases Before and After Low-Impedance Line Preprocessing}
\begin{tabular}{|l|c|c|c|c|}
\hline
\multicolumn{1}{|c|}{\textbf{Case}} & \multicolumn{2}{c|}{\textbf{Before Preprocessing}} & \multicolumn{2}{c|}{\textbf{After Preprocessing}}\\\cline{2-5}
\multicolumn{1}{|c|}{\textbf{Name}} & \textbf{Num.} & \textbf{Num.} & \textbf{Num.} & \textbf{Num.}\\
\multicolumn{1}{|c|}{} & \textbf{\; Buses \;} & \textbf{Lines} & \textbf{\; Buses \;} & \textbf{Lines}\\
\hline
PL-2383wp & 2,383 & 2,869 & 2,177 & 2,690\\\hline
PL-3012wp & 3,012 & 3,572 & 2,292 & 2,851\\\hline
PL-3120sp & 3,120 & 3,693 & 2,314 & 2,886\\\hline
PEGASE-1354 & 1,354 & 1,991 & 1,179 & 1,803\\\hline
PEGASE-2869 & 2,869 & 4,582 & 2,120 & 4,164\\
\hline
\end{tabular}
\label{tab:size}
\end{table}

\begin{figure}[t]
\centering
\includegraphics[totalheight=0.33\textheight]{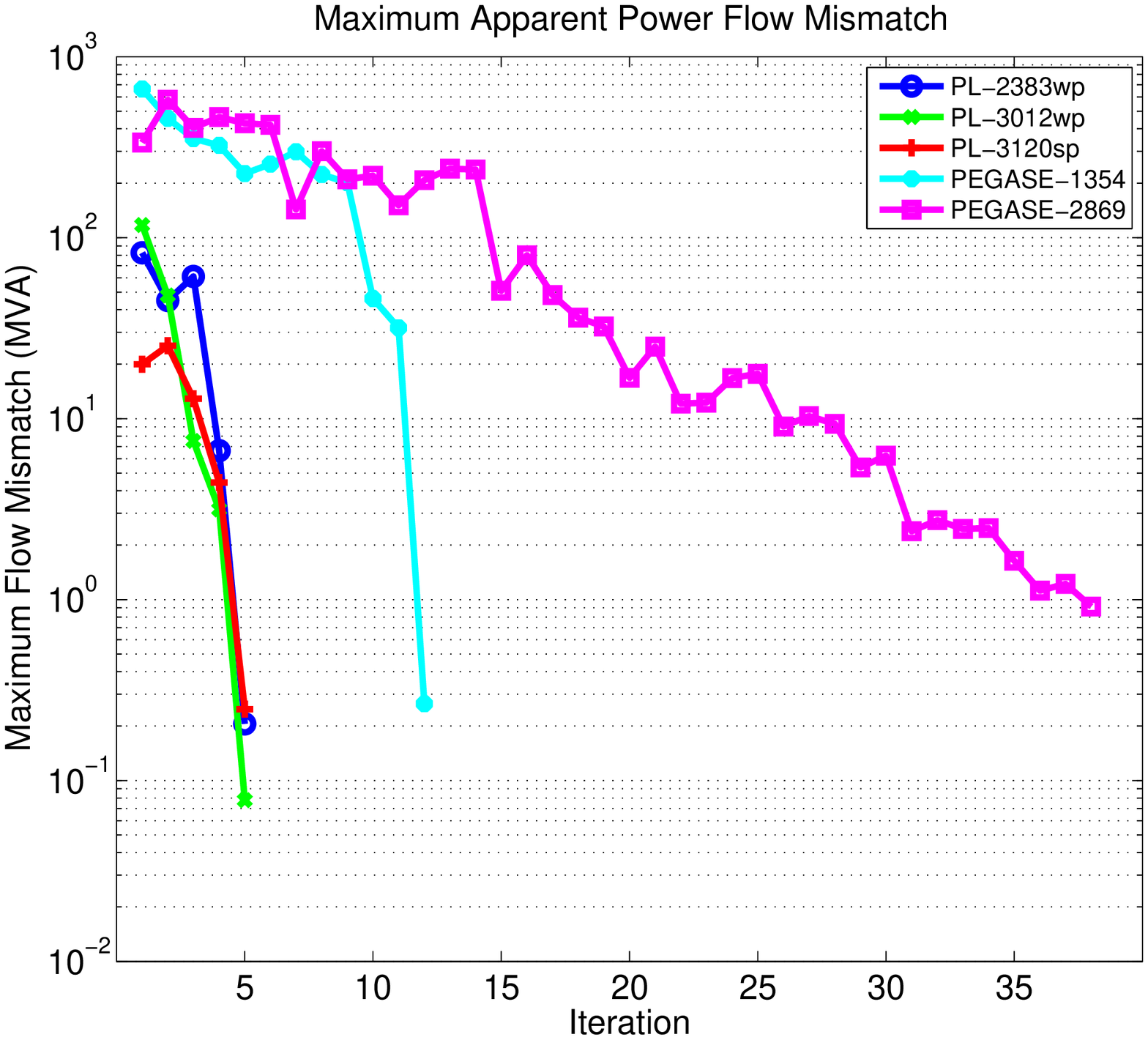}
\caption{Maximum Apparent Power Flow Mismatches versus Iteration of Algorithm~\ref{a:weights} for Active Power Loss Minimizing Test Cases}
\label{f:apparent_flow_mis}
\end{figure}

\begin{figure}[t]
\centering
\includegraphics[totalheight=0.33\textheight]{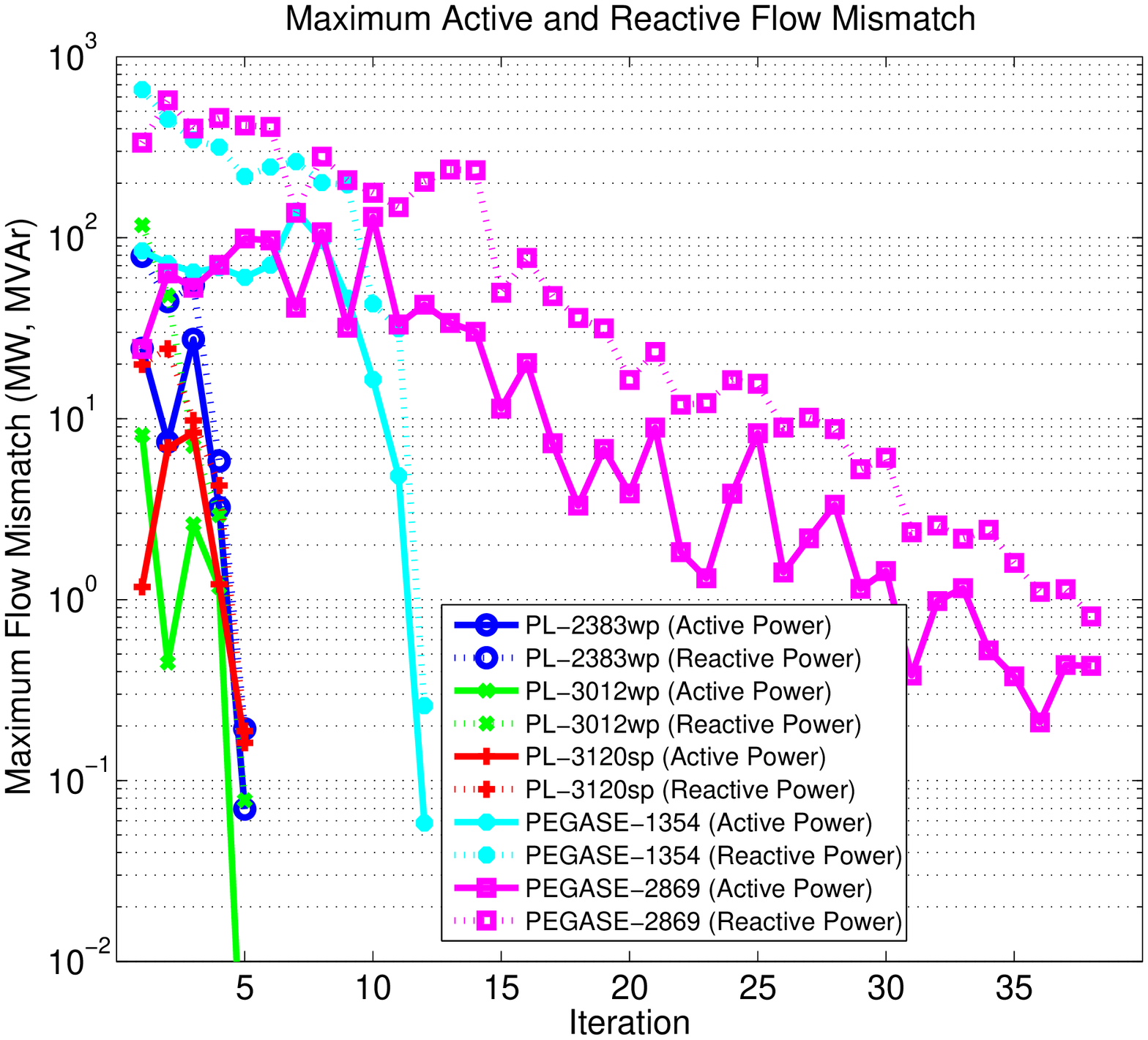}
\caption{Maximum Active and Reactive Power Flow Mismatches versus Iteration of Algorithm~\ref{a:weights} for Active Power Loss Minimizing Test Cases}
\label{f:active_reactive_flow_mis}
\end{figure}

\begin{figure}[t]
\centering
\includegraphics[totalheight=0.33\textheight]{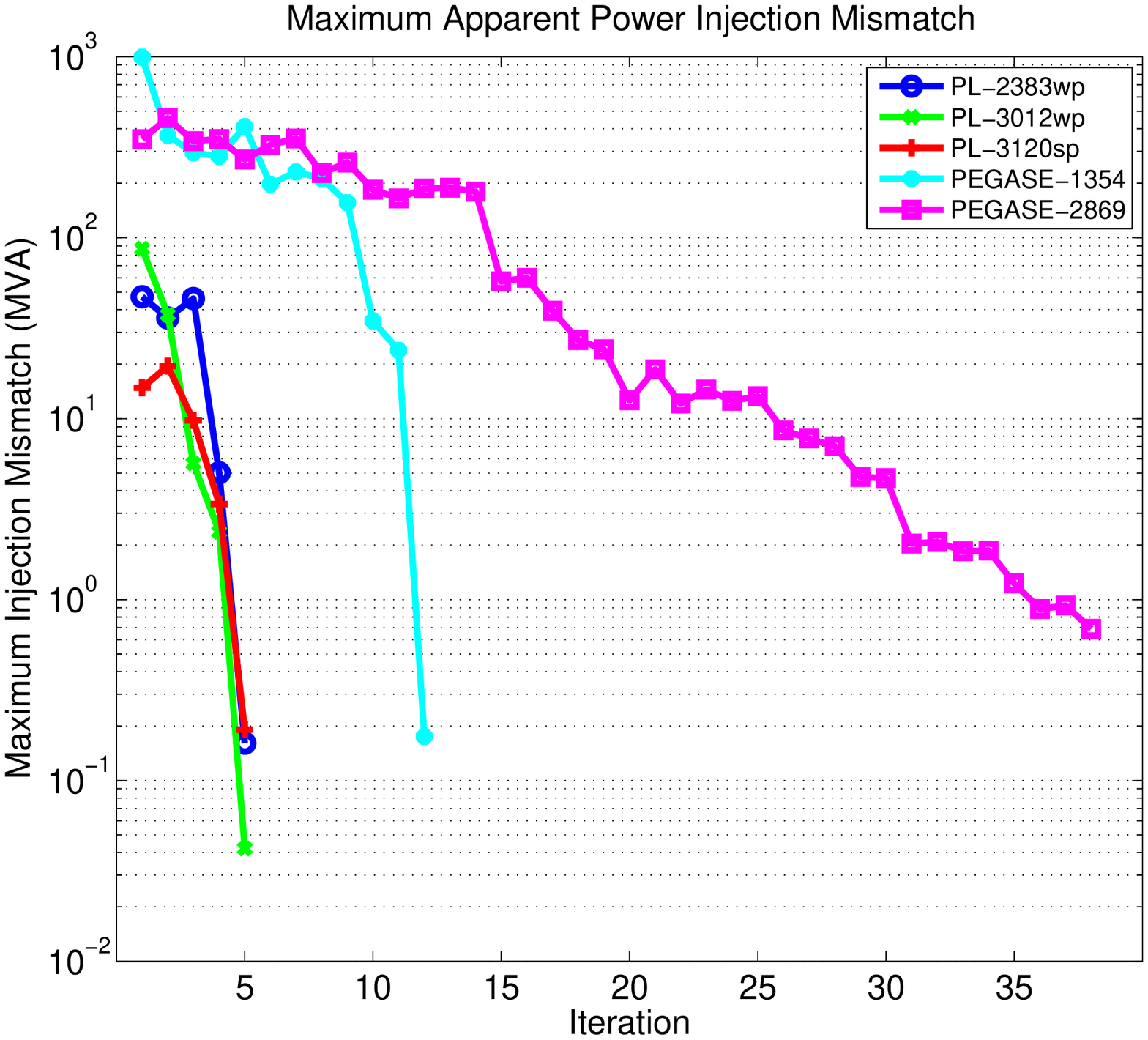}
\caption{Maximum Apparent Power Injection Mismatches versus Iteration of Algorithm~\ref{a:weights} for Active Power Loss Minimizing Test Cases}
\label{f:apparent_injection_mis}
\end{figure}

\begin{figure}[t]
\centering
\includegraphics[totalheight=0.33\textheight]{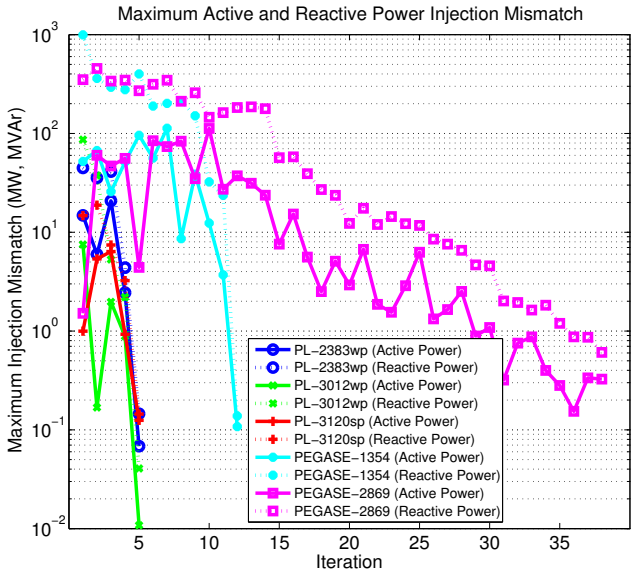}
\caption{Maximum Active and Reactive Power Injection Mismatches versus Iteration of Algorithm~\ref{a:weights} for Active Power Loss Minimizing Test Cases}
\label{f:active_reactive_injection_mis}
\end{figure}

%For these test cases, the number of iterations ranges between 1 and 6. Future work includes altering Algorithm~\ref{a:weights} to more consistently require fewer iterations. 
Note that for the large test cases in Tables~\ref{t:large_gencost_results} and~\ref{t:large_loss_results}, Algorithm~\ref{a:weights} is often computationally faster and has a more straightforward computational implementation than the moment-relaxation-based approaches in~\cite{molzahn_hiskens-sparse_moment_opf,cdc2015}. However, Algorithm~\ref{a:weights} results in feasible points with larger objective values and does not solve as broad a class of OPF problems as existing moment-relaxation-based approaches in~\cite{molzahn_hiskens-sparse_moment_opf,cdc2015}.

% For all test cases, $\delta = 0.5\%$ with the exception of MH118L, which has $\delta = 1\%$ to accommodate a larger relaxation gap of $0.76\%$ for this problem.
Numerical experience suggests that $\delta = 0.5\%$ is
  usually an appropriate parameter choice: as discussed in
  Section~\ref{l:gencost}, the SDP relaxation gap is smaller than
  $0.5\%$ for many test cases. For OPF problems with a significantly
  larger relaxation gap, the proposed approach typically fails to
  yield a feasible solution. Thus, values of $\delta$ that
  differ significantly from $0.5\%$ are not likely to
  be useful in practice.

%The generation cost constraint~\eqref{gencost_constraint_alpha} is binding for all examples in Tables~\ref{t:small_results} and~\ref{t:large_gencost_results}. That is, the feasible points obtained for these examples have objective values that are $\delta$ percent above the lower bound from the SDP relaxation~\eqref{sdpprimal}. Conversely, with the exception of PL-3120sp, the generation cost constraint~\eqref{gencost_constraint_alpha} is not binding for the active power loss minimization results in Table~\ref{t:large_loss_results}. For these test cases, the feasible points obtained from Algorithm~\ref{a:weights} have objective values that range from 0.34\% for \mbox{PL-2383wp} to 0.15\% for \mbox{PEGASE-2869} above the lower bound from~\eqref{sdpprimal}. 

We note that the interior point solver in \matpower{} obtained superior relaxation gaps for the test cases considered in this paper. Within approximately five seconds for the large test cases in Tables~\ref{t:large_gencost_results} and~\ref{t:large_loss_results}, \matpower{} obtained relaxation gaps that ranged from $0.14\%$ to $0.32\%$ smaller than those obtained with Algorithm~\ref{a:weights}. (Of course, \matpower{} cannot provide any measure of the quality of its solution in terms of a lower bound on the globally optimal objective value whereas Algorithm~\ref{a:weights} provides such guarantees.) The smaller relaxation gaps obtained using \matpower{} suggest that smaller values of $\delta$ could be used in Algorithm~\ref{a:weights}. Indeed, additional numerical experiments demonstrated that Algorithm~\ref{a:weights} converged with $\delta = 0.25\%$ (half the value used in previous numerical experiments) for all test cases for which the \matpower{} solution indicated that a value of $\delta = 0.25\%$ was achievable.

We select termination parameter values of
  $\epsilon_{flow}$ and $\epsilon_{inj}$ of 1~MVA, which is a reasonable value for practical power system applications. This tolerance is
  typically numerically achievable with MOSEK's SDP solver, which
  experience suggests is often a limiting factor to obtaining smaller mismatches.
  
  Note that the
maximum mismatches do not necessarily decrease monotonically with each
iteration of Algorithm~\ref{a:weights}. Figs.~\ref{f:apparent_flow_mis}
and~\ref{f:active_reactive_flow_mis} show the maximum flow mismatches (on a logarithmic scale) for the test cases that minimize active power losses (cf.~Table~\ref{t:large_loss_results}). Likewise, Figs.~\ref{f:apparent_injection_mis} and~\ref{f:active_reactive_injection_mis} show the maximum power injection mismatches for the same test cases. Although the mismatches do not always decrease monotonically, there is a generally decreasing trend which results in satisfaction of the termination criteria for each test case. At each iteration, Algorithm~\ref{a:weights} yields larger reactive power mismatches than active power mismatches for these test cases.

%Note that smaller line flow and power injection mismatches can sometimes be obtained by running Algorithm~\ref{a:weights} for several additional iterations. For instance, three extra iterations for MH14L reduces the mismatches below 0.5 MVA.  However, the maximum mismatches do not necessarily decrease monotonically with each iteration of Algorithm~\ref{a:weights}.

Note that it is not straightforward to compare the computational costs of the Laplacian objective approach and other penalization approaches in the literature~\cite{lavaei_mesh,lavaei_allerton2014}. A single solution of the penalized SDP relaxations in~\cite{lavaei_allerton2014} requires approximately the same computational effort as one iteration of Algorithm~\ref{a:weights}. Thus, if one knows appropriate penalty parameters, the method in~\cite{lavaei_allerton2014} is faster for problems where the SDP relaxation is not exact. However, the key advantage of the proposed approach is that there is no need to specify any parameters other than the value of $\delta$ used in the generation cost constraint. In contrast, the literature largely lacks systematic approaches for identifying appropriate parameter values for the penalization methods in~\cite{lavaei_mesh,lavaei_allerton2014}.

\section{Conclusion}
\label{l:conclusion}

The SDP relaxation of~\cite{lavaei_tps} is capable of globally solving
a variety of OPF problems. To address a broader class of OPF problems (i.e., problems for which the SDP relaxation is not exact but provides lower bounds that are close to the global optima),
this paper has described an approach that finds feasible points with
objective values that are within a specified percentage of the global
optimum. Specifically, the approach in this paper adds a constraint to
ensure that the generation cost is within a small specified percentage
of the lower bound obtained from the SDP relaxation. This constraint
frees the objective function to be chosen to yield a \emph{feasible}
(i.e., rank-one) solution rather than a \emph{minimum-cost}
solution. Inspired by previous penalization approaches and results in
the optimization literature, an objective function based on a weighted
Laplacian matrix is selected. The weights for this matrix are
iteratively determined using ``line flow mismatches.'' The proposed
approach is validated through successful application to a variety of
both small and large test cases, including several OPF problems
representing large portions of European power systems. There are,
however, test cases for which the approach takes many iterations to
converge or does not converge at all.

% In contrast to previous work on penalization methods for OPF problems~\cite{lavaei_mesh,lavaei_allerton2014}, this paper presents an algorithmic method for constructing the penalization parameters.

Future work includes modifying the algorithm for choosing the weights
in order to more consistently require fewer iterations. Also, future
work includes testing alternative SDP solution approaches with ``hot
start'' capabilities to improve computational efficiency by leveraging
knowledge of the solution to a ``nearby'' problem from the previous
iteration of the algorithm. Future work also includes extension of the
algorithm to a broader class of OPF problems, such as the test case
WB39mod from~\cite{bukhsh_tps} and several examples in~\cite{kocuk15}.

Additional future work includes leveraging recent results showing that constraints from alternative relaxations (e.g.,~\cite{pscc2014,cedric_msdp,ibm_paper,molzahn_hiskens-sparse_moment_opf,powertech2015,kocuk15,sun2015,coffrin_tightening}) can tighten the SDP relaxation. Augmenting the proposed approach with such constraints may increase its applicability to an broader class of problems.

\section*{Acknowledgment}
The authors acknowledge the support of the Dow Postdoctoral Fellowship in Sustainability, ARPA-E grant \mbox{DE-AR0000232}, and Los Alamos National Laboratory subcontract 270958.

%\ifCLASSOPTIONcaptionsoff
%  \newpage
%\fi

%\IEEEtriggeratref{15}


%\bibliographystyle{IEEEtran}
\bibliographystyle{IEEEtran}
%\bibliography{IEEEabrv,lap_pen}


%\begin{biography}[{\includegraphics[width=1in,height=1.25in,clip,keepaspectratio]{mshell}}]{Michael Shell}

%\begin{IEEEbiography}{Daniel Molzahn}
%\end{IEEEbiography}

%\begin{IEEEbiography}{Bernard Lesieutre}
%\end{IEEEbiography}

%\newpage



%\vfill

\begin{IEEEbiography}[{\includegraphics[width=1in,height=1.25in,clip,keepaspectratio]{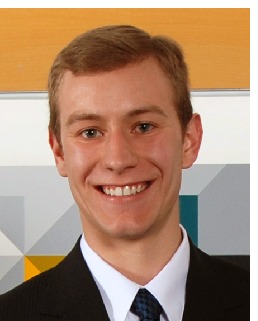}}]{Daniel K. Molzahn}
(S'09-M'13) is a Computational Engineer at Argonne National Laboratory. He recently completed the Dow Postdoctoral Fellow in Sustainability at the University of Michigan, Ann Arbor. He received the B.S., M.S., and Ph.D. degrees in electrical engineering and the Masters of Public Affairs degree from the University of Wisconsin--Madison, where he was a National Science Foundation Graduate Research Fellow. His research focuses on optimization and control of electric power systems.
\end{IEEEbiography}
%\vfill
%\newpage
%

\begin{IEEEbiography}[{\includegraphics[width=1in,height=1.25in,clip,keepaspectratio]{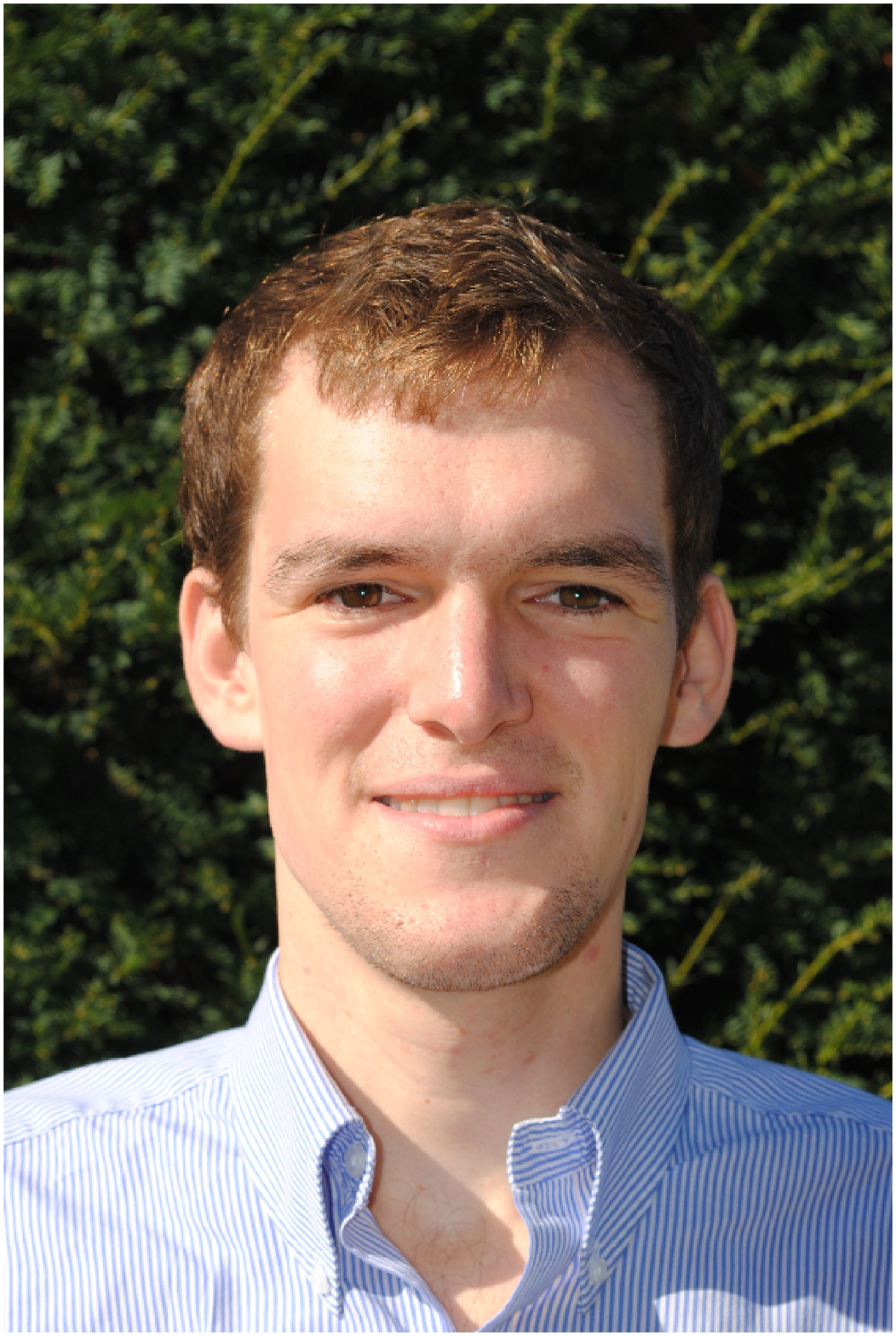}}]{C\'{e}dric Josz}
is pursuing a Ph.D. in Applied Mathematics at Paris VI University in conjunction with RTE, the French transmission system operator, and INRIA, the French national institute in scientific computing. In 2012, he earned a Master's of Engineering degree from ENSTA-Paristech University Paris-Saclay and a Master's degree in optimization from Paris I University.
\end{IEEEbiography}

\begin{IEEEbiography}[{\includegraphics[width=1in,height=1.25in,clip,keepaspectratio]{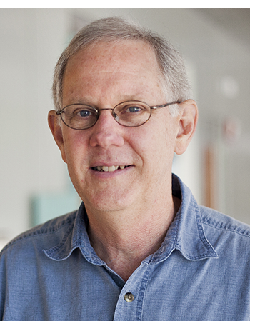}}]{Ian A. Hiskens}
(F'06) is the Vennema Professor of Engineering with the Department of Electrical Engineering and Computer Science, University of Michigan, Ann Arbor. He has held prior appointments in the electricity supply industry (for ten years) and various universities in Australia and the United States. His research focuses on power system analysis, in particular the modelling, dynamics, and control of large-scale, networked, nonlinear systems. His recent activities include integration of renewable generation and new forms of load.

Prof. Hiskens is actively involved in various IEEE societies, and is
VP-Finance of the IEEE System Council. He is a Fellow of Engineers
Australia and a Chartered Professional Engineer in Australia.
\end{IEEEbiography}

\vfill
\newpage
%panciatici_photo
\begin{IEEEbiography}[{\includegraphics[width=1in,height=1.25in,clip,keepaspectratio]{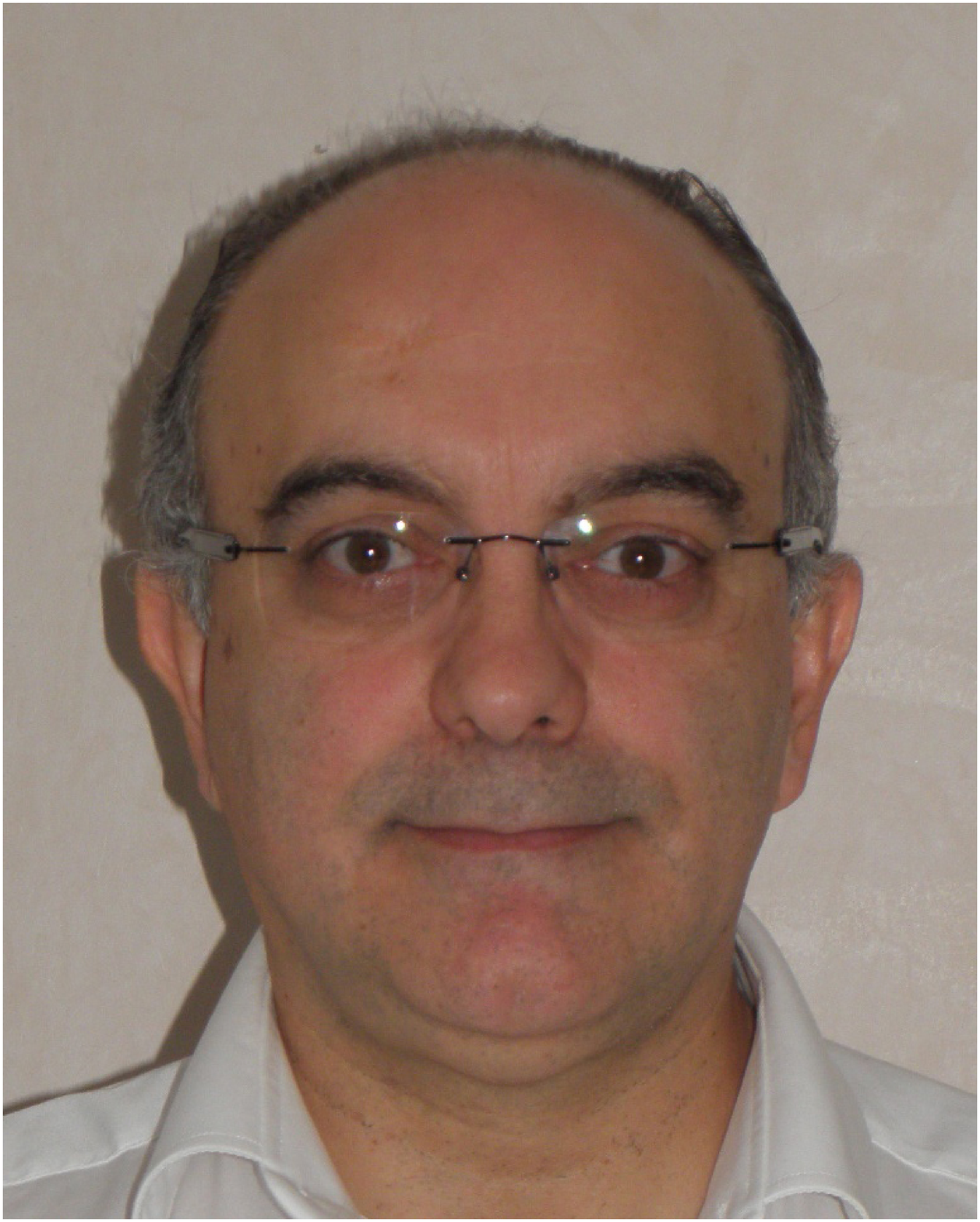}}]{Patrick Panciatici} (M'11) graduated from Supelec, joined EDF R\&D in 1985, and then joined RTE (the French transmission system operator) in 2003. He participated in the creation of the department ``Methods and Support'' at RTE and has more than 25 years of experience in the field of R\&D for transmission systems. Presently, as a scientific advisor, he coordinates and supervises research activities in the Power System Expertise department. Dr. Panciatici is a member of CIGRE, IEEE, and SEE. He participates in the R\&D Plan Working Group of the ENTSO-E (European association of transmission system operators) and in different past and on-going large European Union funded projects (PEGASE, Twenties, iTesla, e-HIGHWAY2050, etc.).
\end{IEEEbiography}
\vfill


\end{document}